\documentclass[12pt]{article}

\usepackage{amsmath}
\usepackage{amssymb}

\newtheorem{thm}{Theorem}[section]

\newtheorem{lem}[thm]{Lemma}

\newcommand{\gp}{\mathfrak{p}}
\newcommand{\gq}{\mathfrak{q}}
\newcommand{\gm}{\mathfrak{m}}

\newcommand{\ga}{\mathfrak{a}}

\newcommand{\zz}{\mathbb{Z}}

\newcommand{\lc}[3]{{\mathrm H}^{#1}_{#2}(#3)} 

\newcommand{\vv}[2]{v_{(#1, #2)}}

\newcommand{\length}[2]{\ell_{#1}(\,#2\,)}
\newcommand{\ass}[2]{{\rm Ass}_{#1}(#2)}

\newcommand{\Min}[2]{{\rm Min}_{#1}(#2)}
\newcommand{\height}[2]{{\rm ht}_{#1}\,#2}
\newcommand{\dep}[2]{{\rm depth}_{#1}\,#2}

\newcommand{\trmatrix}[1]{{}^t\!{#1}} 
\newcommand{\lss}[1]{{}^\ast\!{#1}}
\newcommand{\lsp}[1]{{}^\prime\!{#1}}
\newcommand{\lspp}[1]{{}^{\prime\prime}\!{#1}}
\newcommand{\lran}[1]{\langle{#1}\rangle}  
\newcommand{\lrbr}[1]{[{#1}]} 
\newcommand{\mon}[2]{{\rm m}_{#1}^{#2}} 

\newcommand{\ra}{\longrightarrow}
\newcommand{\homom}[1]{\stackrel{#1}{\longrightarrow}}
\newcommand{\shomom}[1]{\stackrel{#1}{\rightarrow}}

\begin{document}

\title{\Large On a transform of an acyclic complex of length $3$}

\author{
Kosuke Fukumuro, Taro Inagawa and Koji Nishida\footnote{
The last author is supported by KAKENHI (23540042)}
}
\date{
\small
Graduate School of Science, Chiba University, \\
1-33 Yayoi-Cho, Inage-Ku, Chiba-Shi, 263-8522, JAPAN
}
\maketitle

\begin{abstract}
Let $(R\,, \gm)$ be a Cohen-Macaulay local ring and
$Q$ a parameter ideal of $R$\,.
Suppose that an acyclic complex
$F_\bullet$ of length $3$
which is an $R$-free resolution of an ideal $\ga$ of $R$ is given.
In this paper, we describe a concrete procedure to get
an acyclic complex $\lss{F_\bullet}$ of length $3$
that becomes an $R$-free resolution of $\ga :_R Q$\,.
As an application, we compute the symbolic powers of ideals
generated by maximal minors of certain $2 \times 3$ matrices.
\end{abstract}

\section{Introduction}

Let $(R, \gm)$ be a $3$-dimensional Cohen-Macaulay local ring
and $x_1, x_2, x_3$ an sop for $R$\,.
We put $Q = (x_1, x_2, x_3)R$\,.
Suppose that an acyclic complex
\[
0 \ra F_3 \stackrel{\varphi_3}{\ra} F_2 \stackrel{\varphi_2}{\ra} F_1
\stackrel{\varphi_1}{\ra} F_0 = R
\]
of finitely generated free $R$-modules such that
${\rm Im}\,\varphi_3 \subseteq QF_2$ is given.
In this paper, we describe an operation to get an acyclic complex
\[
0 \ra \lss{F_3} \stackrel{\lss{\varphi_3}}{\ra} \lss{F_2} 
\stackrel{\lss{\varphi_2}}{\ra} \lss{F_1} \stackrel{\lss{\varphi_1}}{\ra} \lss{F_0} = R
\]
such that ${\rm Im}\,\lss{\varphi_1} = {\rm Im}\,\varphi_1 :_R Q$ and
${\rm Im}\,\lss{\varphi_3} \subseteq \gm\,\lss{F_2}$\,,
which is called the $\ast$-transform of $F_\bullet$
with respect to $x_1, x_2, x_3$\,.
As we give a practical condition for $\lss{F_3}$ to be vanished,
it is possible to consider when the depth of
$R / ({\rm Im}\,\varphi_1 :_R Q)$ is positive.

If $R$ is a regular local ring and $\ga$ is an ideal of $R$\,,
taking a regular sop and the minimal free resolution of $\ga$
as $x_1, x_2, x_3$ and $F_\bullet$\,, respectively,
we get an acyclic complex $\lss{F_\bullet}$\,,
which is a free resolution of $\ga :_R \gm$\,.
Here, let us notice that we can again take the $\ast$-transform
of $\lss{F_\bullet}$ with respect to
a regular sop since
${\rm Im}\,\lss{\varphi_3} \subseteq \gm\,\lss{F_2}$\,,
and a free resolution of $\ga :_R \gm^2$ is induced.
Repeating this operation, for any $k > 0$\,,
we get a free resolution of $\ga :_R \gm^k$\,.
Moreover, in the case where $\dim R / \ga > 0$\,,
it is possible to find the minimal integer $k > 0$ such that
$\dep{}{R / (\ga :_R \gm^k)} > 0$\,.

As an application of the transformation stated above,
we can compute the symbolic powers of an ideal $I$
generated by the maximal minors of the matrix
\[
\left(
\begin{array}{lll}
x^\alpha & y^\beta & z ^\gamma \\
y^{\beta'} & z^{\gamma'} & x^{\alpha'}
\end{array}
\right)\,,
\]
where $x, y, z$ is an sop for $R$ and
$\alpha, \beta, \gamma, \alpha', \beta', \gamma'$ are positive integers.
As is well known,
$R / I$ is a Cohen-Macaulay local ring with $\dim R / I = 1$\,.
The $n$-th symbolic power of $I$ is defined by
\[
I^{(n)} = \bigcap_{\gp \, \in \, \Min{R}{R / I}}
I^nR_\gp \cap R\,,
\]
and so, if $J$ is an $\gm$-primary ideal such that
$\dep{}{R / (I^n :_R J)} > 0$\,,
we have $I^{(n)} = I^n :_R J$\,.
Because $I$ is generated by a d-sequence,
it is possible to describe a minimal free resolution of $I^n$\,,
and we can take its
$\ast$-transform with respect to
$x^{\alpha''}, y^{\beta''}, z^{\gamma''}$\,,
where $\alpha'' = \min\{\,\alpha, \alpha'\,\},
\beta'' = \min\{\,\beta, \beta'\,\}$ and $\gamma'' = \min\{\,\gamma, \gamma'\,\}$\,.
If the length of the resulting acyclic complex is still $3$\,,
we have $\dep{}{R / (I^n :_R (x^{\alpha''}, y^{\beta''}, z^{\gamma''}))} = 0$\,,
and it means $I^n :_R (x^{\alpha''}, y^{\beta''}, z^{\gamma''}) \subsetneq I^{(n)}$\,.
Then, we take once more the $\ast$-transform with respect to
suitable powers of $x, y, z$\,.
By repeating such operation several times,
we can reach $I^{(n)}$\,.
In Section 4, we carry it out for $n = 2, 3$\,.
Our method is completely different from that used in 
\cite{gns} and \cite{ksz},
in which $I^{(2)}$ and $I^{(3)}$ are computed in the case where
$R$ is a regular local ring with the maximal ideal $\gm = (x, y, z)R$ and
$I$ is a prime ideal.

In the last section,
assuming that $R$ is a regular local ring and
$x, y, z$ is a regular sop for $R$\,,
we compute the length of $I^{(n)} / I^n$
for all $n \geq 1$ in the case where
$I$ is generated by the maximal minors of the matrix
\[
\left(\begin{array}{lll}
x & y & z \\
y & z & x^2
\end{array}\right)\,.
\]
Starting from a minimal free resolution of $I^n$\,,
we repeat the operation to take the $\ast$-transform
with respect to $x, y, z$ successively.
Then we get a free resolution of $I^n :_R \gm^k$
for any $k \geq 1$\,,
and it follows that
$I^{(n)} = I^n :_R \gm^q$\,,
where $q$ is the largest integer with
$q \leq n/2$\,.
As our method enables us to compute the length of
$(I^n :_R \gm^k) / (I^n :_R \gm^{k-1})$
for any $k \geq 1$\,,
we get to know the length of $I^{(n)} / I^n$ exactly\,.
As a consequence of our result, we see
that the $\epsilon$-multiplicity of $I$\,,
which is an invariant defined by
\[
\epsilon(I) := \lim_{n \rightarrow \infty}
\frac{3!}{n^3}\!\cdot\!\length{R}{I^{(n)} / I^n}
\]
(cf. \cite{chs}, \cite{uv}), coincides with $1/2$\,.

Throughout this paper,
$R$ is a $3$-dimensional Cohen-Macaulay local ring
with the maximal ideal $\gm$\,.
If an $R$-module is a direct sum of finite number of $R$-modules,
its elements are denoted by column vectors.
In particular, if an $R$-module $F$ is a direct sum of
two $R$-modules, say $F = G \oplus H$\,,
the elements in $F$ of the form
\[
\binom{g}{0} \hspace{1ex}\mbox{and}\hspace{1ex} \binom{0}{h}
\hspace{3ex}
( g \in G,\, h \in H )
\]
are denoted by $\lrbr{g}$ and $\lran{h}$\,, respectively.
Moreover, if $V$ is a subset of $G$\,,
then the family $\{\, \lrbr{g} \,\}_{g \in V}$
is denoted by $\lrbr{V}$\,.
Similarly $\lran{W}$ is defined for a subset $W$ of $H$\,.
If $T$ is a subset of an $R$-module,
we denote by $R\!\cdot\! T$ the $R$-submodule generated by $T$\,.

\section{Preliminaries}

In this section, we summarize preliminary results.
Although they might be well known,
we give the proofs for completeness.
Let us begin with the following lemma.

\begin{lem}\label{2.1}
Let $G_\bullet$ and $F_\bullet$ be acyclic complexes,
whose boundary maps are denoted by
$\partial_\bullet$ and $\varphi_\bullet$\,, respectively.
Suppose that a chain map $\sigma_\bullet : G_\bullet \ra F_\bullet$
is given and $\sigma_0^{-1}({\rm Im}\,\varphi_1) =
{\rm Im}\,\partial_1$ holds.
Then the mapping cone ${\rm Con}(\sigma_\bullet)$ is acyclic.
Hence, if $G_\bullet$ and $F_\bullet$ are complexes of
finitely generated free $R$-modules,
then ${\rm Con}(\sigma_\bullet)$ gives an $R$-free resolution of
${\rm Im}\,\varphi_1 + {\rm Im}\,\sigma_0$\,.
\end{lem}

\noindent
{\it Proof}.\,
We put $C_\bullet = {\rm Con}(\sigma_\bullet)$
and denote its boundary map by $d_\bullet$\,. Then
\[
C_i = F_i \oplus G_{i-1} \hspace{1ex}, \hspace{3ex}
d_i =
\left(\begin{array}{cc}
\varphi_i & (-1)^{i-1}\sigma_{i-1} \\
0 & \partial_{i-1}
\end{array}\right)
\]
for any $i \geq 2$ and
\[
C_1 = F_1 \oplus G_0 \hspace{1ex}, \hspace{3ex}
C_0 = F_0 \hspace{1ex}, \hspace{3ex}
d_1 = 
\left(\begin{array}{cc}
\varphi_1 & \sigma_0
\end{array}\right)\,.
\]
Let us take any
\[
\binom{x_1}{y_0} \in {\rm Ker}\,d_1\,.
\]
Then $\varphi_1(x_1) + \sigma_0(y_0) = 0$\,,
so $y_0 \in \sigma_0^{-1}({\rm Im}\,\varphi_1) = {\rm Im}\,\partial_1$\,.
Hence there exists $y_1 \in G_1$ such that $y_0 = \partial_1(y_1)$\,.
We have $\varphi_1(\partial_1(y_1)) = \sigma_0(\partial_1(y_1)) =
\sigma_0(y_0) = -\varphi_1(x_1)$\,, and so,
$\varphi_1(x_1 + \sigma_1(y_1)) = 0$\,.
This means that there exists $x_2 \in F_2$ such that
$x_1 + \sigma_1(y_1) = \varphi_2(x_2)$ as $\lc{}{1}{F_\bullet} = 0$\,.
Then
\[
\binom{x_1}{y_0} = \binom{\varphi_2(x_2) - \sigma_1(y_1)}{\partial_1(y_1)} =
d_2 \binom{x_2}{y_1} \in {\rm Im}\,d_2\,.
\]
Therefore we get $\lc{}{1}{C_\bullet} = 0$\,.
On the other hand, the exact sequence
\[
0 \ra G_\bullet(-1) \ra C_\bullet \ra F_\bullet \ra 0
\]
of complexes induces a long exact sequence
\[
\cdots \ra \lc{}{i-1}{G_\bullet} \ra \lc{}{i}{C_\bullet}
\ra \lc{}{i}{F_\bullet} \ra \cdots
\]
of homology modules.
If $i \geq 2$\,, $\lc{}{i-1}{G_\bullet} = \lc{}{i}{F_\bullet} = 0$\,,
and so, $\lc{}{i}{C_\bullet} = 0$\,.
Thus we have seen that $C_\bullet$ is an acyclic complex.

\begin{lem}\label{2.2}
Let $C_{\bullet\bullet}$ be a double complex
such that $C_{pq} = 0$ unless $0 \leq p, q \leq 3$\,.
For any $p, q \in \zz$\,,
we denote the boundary maps $C_{pq} \ra C_{p-1, q}$
and $C_{pq} \ra C_{p, q-1}$ by $d'_{pq}$ and $d''_{pq}$\,, respectively.
If $C_{p\bullet}$ and $C_{\bullet q}$
are acyclic for any $0 \leq p, q \leq 3$\,,
we have the following assertions on the total complex $T_\bullet$\,,
whose boundary map is denoted by $d_\bullet$\,.
\begin{itemize}
\item[{\rm (1)}]
Suppose that $\xi_0 \in C_{0, 3}$ and $\xi_1 \in C_{1, 2}$
such that $d''_{0, 3}(\xi_0) + d'_{1, 2}(\xi_1) = 0$ are given.
Then there exist $\xi_2 \in C_{2, 1}$ and $\xi_3 \in C_{3, 0}$
such that
$\trmatrix{
\begin{pmatrix} \xi_0 & \xi_1 & \xi_2 & \xi_3 \end{pmatrix}}
\in {\rm Ker}\,d_3
\subseteq T_3 = C_{0, 3} \oplus C_{1, 2} \oplus C_{2, 1} \oplus C_{3, 0}$\,.
\item[{\rm (2)}]
Let $\trmatrix{
\begin{pmatrix} \xi_0 & \xi_1 & \xi_2 & \xi_3 \end{pmatrix}}
\in {\rm Ker}\,d_3 \subseteq T_3
= C_{0, 3} \oplus C_{1, 2} \oplus C_{2, 1} \oplus C_{3, 0}$
and let $\xi_3 \in {\rm Im}\,d''_{3, 1}$\,.
Then $\trmatrix{
\begin{pmatrix} \xi_0 & \xi_1 & \xi_2 & \xi_3 \end{pmatrix}}
\in {\rm Im}\,d_4$\,.
In particular, $\xi_0 \in {\rm Im}\,d'_{1, 3}$\,.
\end{itemize}
\end{lem}

\noindent
{\it Proof}. \,
(1) \,
Because $d'_{1,1}(d''_{1, 2}(\xi_1)) = d''_{0, 2}(d'_{1, 2}(\xi_1)) =
d''_{0, 2}(-d''_{0, 3}(\xi_0)) = 0$\,,
it follows that
$d''_{1, 2}(\xi_1) \in {\rm Ker}\,d'_{1, 1} = {\rm Im}\,d'_{2, 1}$\,.
So, we can write $d''_{1, 2}(\xi_1) = d'_{2, 1}(\xi_2)$\,,
where $\xi_2 \in C_{2, 1}$\,.
Then $d'_{2, 0}(d''_{2, 1}(\xi_2)) = d''_{1, 1}(d'_{2, 1}(\xi_2)) =
d''_{1, 1}(d''_{1, 2}(\xi_1)) = 0$\,,
and so $d''_{2, 1}(\xi_2) \in {\rm Ker}\,d'_{2, 0} = {\rm Im}\,d'_{3, 0}$\,.
Hence there exists $\xi_3 \in C_{3, 0}$ such that
$d''_{2, 1}(\xi_2) = -d'_{3, 0}(\xi_3)$\,.
Now we have
\[
d_3
\begin{pmatrix}
\xi_0 \\
\xi_1 \\
\xi_2 \\
\xi_3
\end{pmatrix} = 
\begin{pmatrix}
d''_{0, 3}(\xi_0) + d'_{1, 2}(\xi_1) \\
-d''_{1, 2}(\xi_1) + d'_{2, 1}(\xi_2) \\
d''_{2, 1}(\xi_2) + d'_{3, 0}(\xi_3))
\end{pmatrix} = 
\begin{pmatrix}
0 \\
0 \\
0
\end{pmatrix}\,,
\]
and so we get the required assertion.

(2) \,
Because $\trmatrix{
\begin{pmatrix} \xi_0 & \xi_1 & \xi_2 & \xi_3 \end{pmatrix}}
\in {\rm Ker}\,d_3$\,, we have
(i)\, $d''_{0, 3}(\xi_0) + d'_{1, 2}(\xi_1) = 0$\,,
(ii)\, $-d''_{1, 2}(\xi_1) + d'_{2, 1}(\xi_2) = 0$\, and
(iii)\, $d''_{2, 1}(\xi_2) + d'_{3, 0}(\xi_3) = 0$\,.
On the other hand, as $\xi_3 \in {\rm Im}\,d''_{3, 1}$\,,
there exists $\eta_3 \in C_{3, 1}$ such that
$\xi_3 = -d''_{3, 1}(\eta_3)$\,.
Then, as $d'_{3, 0}(\xi_3) = -d'_{3, 0}(d''_{3, 1}(\eta_3)) =
-d''_{2, 1}(d'_{3, 1}(\eta_3))$\,,
by (iii) we get
$\xi_2 - d'_{3, 1}(\eta_3) \in {\rm Ker}\,d''_{2, 1} = {\rm Im}\,d''_{2, 2}$\,.
So, we can write
$\xi_2 - d'_{3, 1}(\eta_3) = d''_{2, 2}(\eta_2)$\,,
where $\eta_2 \in C_{2, 2}$\,.
Then $d'_{2, 1}(\xi_2) =
d'_{2, 1}(d'_{3, 1}(\eta_3) + d''_{2, 2}(\eta_2)) =
d'_{2, 1}(d''_{2, 2}(\eta_2)) = d''_{1, 2}(d'_{2, 2}(\eta_2))$\,,
and so by (ii) we get $\xi_1 - d'_{2, 2}(\eta_2) \in
{\rm Ker}\,d''_{1, 2} = {\rm Im}\,d''_{1, 3}$\,.
Hence there exists $\eta_1 \in C_{1, 3}$ such that
$\xi_1 - d'_{2, 2}(\eta_2) = -d''_{1, 3}(\eta_1)$\,.
This means $d'_{1, 2}(\xi_1) =
d'_{1, 2}(d'_{2, 2}(\eta_2) - d''_{1, 3}(\eta_1)) =
-d'_{1, 2}(d''_{1, 3}(\eta_1)) = -d''_{0, 3}(d'_{1, 3}(\eta_1))$\,,
and so by (i) we get
$\xi_0 - d'_{1, 3}(\eta_1) \in {\rm Ker}\,d''_{0, 3}$\,.
However, as $C_{0, \bullet}$ is acyclic and $C_{0, 4} = 0$\,,
$d''_{0, 3}$ is injective.
Hence $\xi_0 = d'_{1, 3}(\eta_1)$\,.
Now we have
\[
d_4
\begin{pmatrix}
\eta_1 \\
\eta_2 \\
\eta_3
\end{pmatrix} = 
\begin{pmatrix}  
d'_{1, 3}(\eta_1) \\
-d''_{1, 3}(\eta_1) + d'_{2, 2}(\eta_2) \\
d''_{2, 2}(\eta_2) + d'_{3, 1}(\eta_3) \\
-d''_{3, 1}(\eta_3)
\end{pmatrix} =
\begin{pmatrix}
\xi_0 \\
\xi_1 \\
\xi_2 \\
\xi_3
\end{pmatrix}\,,
\]
and the proof is completed.

\begin{lem}\label{2.3}
Suppose that
\[
0 \ra F \stackrel{\varphi}{\ra} G \stackrel{\psi}{\ra} H \stackrel{\rho}{\ra} L
\]
is an exact sequence of $R$-modules.
Then the following assertions hold.

\begin{itemize}

\item[{\rm (1)}]
If there exists a homomorphism $\phi : G \ra F$ of $R$-modules
such that $\phi \circ \varphi = {\rm id}_F$\,, then
\[
0 \ra \lss{G} \stackrel{\lss{\psi}}{\ra} H \stackrel{\rho}{\ra} L
\]
is exact, where $\lss{G} = {\rm Ker}\,\phi$ and
$\lss{\psi}$ is the restriction of $\psi$ to $\lss{G}$ .

\item[{\rm (2)}]
If $F = \lsp{F} \oplus \lss{F}$\,, $G = \lsp{G} \oplus \lss{G}$\,,
$\varphi (\lsp{F}) = \lsp{G}$ and $\varphi (\lss{F}) \subseteq \lss{G}$\,, then
\[
0 \ra \lss{F} \stackrel{\lss{\varphi}}{\ra} \lss{G} \stackrel{\lss{\psi}}{\ra} H
\stackrel{\rho}{\ra} L
\]
is exact, where $\lss{\varphi}$ and $\lss{\psi}$ are
the restrictions of $\varphi$ and $\psi$
to $\lss{F}$ and $\lss{G}$\,, respectively.

\end{itemize}

\end{lem}

\noindent
{\it Proof}. \,
(1) \,
First, we take any $u \in {\rm Ker}\,\rho$\,.
As ${\rm Im}\,\psi = {\rm Ker}\,\rho$\,,
there exists $v \in G$ such that $\psi(v) = u$\,.
Then, setting $\lss{v} = v - \varphi(\phi(v))$\,,
we have $\lss{v} \in \lss{G}$ and
$\psi(\lss{v}) = \psi(v) = u$\,, and so $u \in \psi(\lss{G})$\,.
Hence ${\rm Ker}\,\rho = {\rm Im}\,\lss{\psi}$\,.

Next, we take any $\lss{v} \in \lss{G}$ such that
$\psi(\lss{v}) = 0$\,.
As ${\rm Im}\,\varphi = {\rm Ker}\,\psi$\,,
there exists $w \in F$ such that $\varphi(w) = \lss{v}$\,.
Then $w = \phi(\varphi(w)) = \phi(\lss{v}) = 0$\,,
and so $\lss{v} = \varphi(0) = 0$\,.
Hence $\lss{\varphi}$ is injective.

(2) \,
First, we take any $u \in {\rm Ker}\,\rho$\,.
As ${\rm Im}\,\psi = {\rm Ker}\,\rho$\,,
there exists $v \in G$ such that $\psi(v) = u$\,.
We write $v = \lsp{v} + \lss{v}$\,, where $\lsp{v} \in \lsp{G}$
and $\lss{v} \in \lss{G}$\,,
and choose $\lsp{w} \in \lsp{F}$ so that $\varphi(\lsp{w}) = \lsp{v}$\,.
Then we have $u = \psi(\varphi(\lsp{w}) + \lss{v}) =
\psi(\lss{v}) \in \psi(\lss{G})$\,.
Hence ${\rm Im}\,\lss{\psi} = {\rm Ker}\,\rho$\,.

Next, we take any $\lss{v} \in \lss{G}$ such that
$\psi(\lss{v}) = 0$\,.
As ${\rm Im}\,\varphi = {\rm Ker}\,\psi$\,,
there exists $w \in F$ such that
$\varphi(w) = \lss{v}$\,.
We write $w = \lsp{w} + \lss{w}$\,,
where $\lsp{w} \in \lsp{F}$ and $\lss{w} \in \lss{F}$\,.
Then $\lss{v} - \varphi(\lss{w}) =
\varphi(\lsp{w}) \in \lsp{G} \cap \lss{G} = 0$\,.
This means $\lss{v} = \varphi(\lss{w}) \in \varphi(\lss{G})$\,.
Hence ${\rm Im}\,\lss{\varphi} = {\rm Ker}\,\lss{\psi}$\,,
and the proof is completed.

\section{$\ast$-transform}

Let $R$ be a $3$-dimensional Cohen-Macaulay local ring
with the maximal ideal $\gm$ and
\[
0 \homom{} F_3 \homom{\varphi_3} F_2 \homom{\varphi_2}
F_1 \homom{\varphi_1} F_0 = R
\]
be an acyclic complex of finitely generated free $R$-modules
such that ${\rm Im}\,\varphi_3 \subseteq QF_2$\,,
where $Q = (x_1, x_2, x_3)$ is a parameter ideal of $R$\,.
We put $\ga = {\rm Im}\,\varphi_1$\,,
which is an ideal of $R$\,.
In this section, transforming $F_\bullet$ suitably,
we aim to construct an acyclic complex
\[
0 \homom{} \lss{F_3} \homom{\lss{\varphi_3}}
\lss{F_2} \homom{\lss{\varphi_2}}
\lss{F_1} \homom{\lss{\varphi_1}} \lss{F_0} = R
\]
of finitely generated free $R$-modules such that
${\rm Im}\,\lss{\varphi_3} \subseteq \gm\lss{F_2}$ and
${\rm Im}\,\lss{\varphi_1} = \ga :_R Q$\,.
Let us call $\lss{F_\bullet}$ the $\ast$-transform of $F_\bullet$
with respect to $x_1, x_2, x_3$\,.

In this operation, we use the Koszul complex
$K_\bullet = K_\bullet(x_1, x_2, x_3)$\,.
Let $e_1, e_2, e_3$ be an $R$-free basis of $K_1$ and set
$\check{e}_1 = e_2 \wedge e_3$\,,
$\check{e}_2 = e_1 \wedge e_3$\,, $\check{e}_3 = e_1 \wedge e_2$\,.
Then $\check{e}_1, \check{e}_2, \check{e}_3$ is an $R$-free basis of $K_2$\,. 
Furthermore $e_1 \wedge e_2 \wedge e_3$ is an $R$-free basis of $K_3$\,.  
The boundary maps of $K_\bullet$
which we denote by $\partial_\bullet$ satisfy
\begin{eqnarray*}
\partial_1(e_i) & = & x_i \hspace{2ex} \mbox{for any} \hspace{1ex} i = 1, 2, 3\,,  \\
\partial_2(e_i \wedge e_j) & = & 
x_i e_j - x_j e_i \hspace{2ex} \mbox{if} \hspace{1ex} 1 \leq i < j \leq 3\,, \\
\partial_3(e_1 \wedge e_2 \wedge e_3) &  = & 
x_1 \check{e}_1 - x_2 \check{e}_2 + x_3 \check{e}_3\,.
\end{eqnarray*}
As $x_1, x_2, x_3$ is an $R$-regular sequence,
$K_\bullet$ gives an $R$-free resolution of $R / Q$\,.
Hence, for any $R$-module $M$\,,
We have the following commutative diagram;
\[
\begin{array}{cccccccl}
{\rm Hom}_R(K_2\,, M) & \homom{{\rm Hom}_R(\partial_3\,,\, M)} &
   {\rm Hom}_R(K_3\,, M) & \ra & {\rm Ext}_R^3(R / Q\,, M) & \ra & 0 & {\rm (ex)} \\
\downarrow\vcenter{\rlap{$\scriptstyle{\cong}$}} &   &
  \downarrow\vcenter{\rlap{$\scriptstyle{\cong}$}} &   &   &   &    &    \\
M^{\oplus 3} & \homom{(x_1 \hspace{1ex} -x_2 \hspace{1ex} x_3)} & 
    M & \ra & M / QM & \ra & 0 & {\rm (ex)\,,}
\end{array}
\]
which implies ${\rm Ext}_R^3(R / Q\,, M) \cong M / QM$\,.
Using this fact, we show the next result.

\begin{thm}\label{3.1}
$(\ga :_R Q) / \ga \cong F_3 / QF_3$\,.
\end{thm}

\noindent
{\it Proof}.\,
As is well known, we have
\[
(\ga :_R Q) / \ga \cong {\rm Hom}_R(R / Q\,, R / \ga)\,.
\]
Applying ${\rm Hom}_R(R / Q, \,\bullet\,)$ to the exact sequence
$0 \ra \ga \ra R \ra R / \ga \ra 0$\,,
we get
\[
{\rm Hom}_R(R / Q\,, R / \ga) \cong {\rm Ext}_R^1(R / Q\,, \ga)
\]
as ${\rm Hom}_R(R / Q\,, R) = {\rm Ext}_R^1(R / Q\,, R) = 0$\,.
We put $L = {\rm Ker}\,\varphi_1 = {\rm Im}\,\varphi_2$
and consider the exact sequences
\[
0 \ra L \ra F_1 \homom{\varphi_1} \ga \ra 0 \hspace{2ex}
\mbox{and} \hspace{2ex}
0 \ra F_3 \homom{\varphi_3} F_2 \ra L \ra 0\,.
\]
Because ${\rm Ext}_R^i(R / Q\,, F_j) = 0$ for any
$0 \leq i \leq 2$ and $0 \leq j \leq 3$\,, we get
\[
{\rm Ext}_R^1(R / Q\,, \ga) \cong {\rm Ext}_R^2(R / Q\,, L)
\]
and the commutative diagram
\[
\begin{array}{cccccccl}
0 & \ra & {\rm Ext}_R^2(R / Q\,, L) & \ra & {\rm Ext}_R^3(R / Q\,, F_3) &
  \homom{\tilde{\varphi}_3} & {\rm Ext}_R^3(R / Q\,, F_2) & {\rm (ex)} \\
  &   &   &   & 
  \downarrow\vcenter{\rlap{$\scriptstyle{\cong}$}} &    &
    \downarrow\vcenter{\rlap{$\scriptstyle{\cong}$}}   \\
  &   &   &   & F_3 / QF_3 & \homom{\overline{\varphi}_3} & 
     F_2 / QF_2 & {\rm (ex)\,,}
\end{array}
\]
where $\tilde{\varphi}_3$ and $\overline{\varphi}_3$
denote the maps induced from $\varphi_3$\,.
Let us notice $\overline{\varphi}_3 = 0$ as
${\rm Im}\,\varphi_3 \subseteq QF_2$\,.
Hence
\[
{\rm Ext}_R^2(R / Q\,, L) \cong F_3 / QF_3\,,
\]
and so the required isomorphism follows.

\vspace{1.5em}

Let us fix an $R$-free basis of $F_3$\,, say
$\{\,w_\lambda\,\}_{\lambda \in \Lambda}$\,.
Let $\{\,v_{(i, \lambda)}\,\}_{1\leq i \leq 3\,,\, \lambda \in \Lambda}$
be a family of elements in $F_2$ such that
\[
\varphi_3(w_\lambda) = \sum_{i = 1}^3 x_i \!\cdot\! v_{(i, \lambda)}
\]
for any $\lambda \in \Lambda$\,.
We put
$\tilde{\Lambda} = \{\,1\,, 2\,, 3\,\} \times \Lambda$\,.
The next result is the essential part of the process to get
$\lss{F_\bullet}$\,.

\begin{thm}\label{3.2}
There exists a chain map 
$\sigma_\bullet : K_\bullet \otimes_R F_3 \ra F_\bullet$\,
\[
{
\begin{array}{ccccccccc}
0 & \homom{} & K_3 \otimes_R F_3 & \homom{\partial_3 \otimes {\rm id}} &
   K_2 \otimes_R F_3 & \homom{\partial_2 \otimes {\rm id}} & K_1 \otimes_R F_3 &
        \homom{\partial_1 \otimes {\rm id}} & K_0 \otimes_R F_3 \\
   &   & 
\downarrow\vcenter{%
 \rlap{$\scriptstyle{\sigma_3}$}}
 &    & 
\downarrow\vcenter{%
 \rlap{$\scriptstyle{\sigma_2}$}}
 &   & 
\downarrow\vcenter{%
 \rlap{$\scriptstyle{\sigma_1}$}}
 &   & 
\downarrow\vcenter{%
 \rlap{$\scriptstyle{\sigma_0}$}}
  \\
0 & \homom{} & F_3 & \homom{\varphi_3} & F_2 & \homom{\varphi_2} & 
   F_1 & \homom{\varphi_1} & F_0
\end{array}
}
\]
satisfying the following conditions. 
\begin{itemize}
\item[{\rm (1)}]
$\sigma_0^{-1}({\rm Im}\,\varphi_1) = {\rm Im}\,(\partial_1 \otimes {\rm id}_{F_3})$\,.  
\item[{\rm (2)}]
${\rm Im}\,\sigma_0 + {\rm Im}\,\varphi_1 = \ga :_R Q$\,. 
\item[{\rm (3)}]
$\sigma_2(\check{e}_i \otimes w_\lambda) = (-1)^i\!\cdot\!\vv{i}{\lambda}$
for any $(i, \lambda) \in \tilde{\Lambda}$\,. 
\item[{\rm (4)}]
$\sigma_3((e_1 \wedge e_2 \wedge e_3) \otimes w_\lambda) =
-w_\lambda$ for any $\lambda \in \Lambda$\,.
\end{itemize}
\end{thm}

\noindent
{\it Proof}.\,
Let us consider the double complex $K_\bullet \otimes_R F_\bullet$\,.
\[
\begin{array}{ccccccc}
  &  & \vdots &  & \vdots &  &  \\
  &  & \big\downarrow &  & \big\downarrow &  &  \\
\cdots & \ra & K_i \otimes_R F_j & 
   \homom{\partial_i\,\otimes\,{\rm id}_{F_j}} & K_{i-1} \otimes_R F_j & \ra & \cdots \\
  &  &  &  &  &  &  \\
  &  & \big\downarrow\vcenter{\rlap{$\scriptstyle{{\rm id}_{K_i}\,\otimes\,\varphi_j}$}} &
    & \big\downarrow\vcenter{\rlap{$\scriptstyle{{\rm id}_{K_{i-1}}\,\otimes\,\varphi_j}$}} &
      &    \\
\cdots & \ra & K_i \otimes_R F_{j-1} & 
   \homom{\partial_i\,\otimes\,{\rm id}_{F_{j-1}}} & K_{i-1} \otimes_R F_{j-1} & \ra & \cdots \\
  &  &  &  &  &  &  \\
  &  & \big\downarrow &   & \big\downarrow &  &   \\
  &  & \vdots &    & \vdots &  &  
\end{array}
\]
We can take it as $C_{\bullet\bullet}$ in \ref{2.2}.
Let $T_\bullet$ be the total complex and $d_\bullet$ be its boundary map.
For any $\lambda \in \Lambda$\,, we set
$\xi_0(\lambda) = 1 \otimes w_\lambda \in K_0 \otimes F_3$ and
$\xi_1(\lambda) = 
-\sum_{i = 1}^3\,e_i \otimes v_{(i, \lambda)} \in K_1 \otimes_R F_2$\,.
Then we have
\[
({\rm id}_{K_0} \otimes \varphi_3)(\xi_0(\lambda)) +
  (\partial_1 \otimes {\rm id}_{F_2})(\xi_1(\lambda)) = 0\,,
\]
and so by \ref{2.2} there exist $\xi_2(\lambda) \in K_2 \otimes_R F_1$
and $\xi_3(\lambda) \in K_3 \otimes_R F_0$ such that
\[
\trmatrix{\,(\,\xi_0(\lambda) \hspace{1ex} \xi_1(\lambda) \hspace{1ex}
  \xi_2(\lambda) \hspace{1ex} \xi_3(\lambda)\,)} \in
{\rm Ker}\,(T_3 \homom{d_3} T_2)\,.
\]
When this is the case, the following equalities hold;
\begin{align*}
\mbox{(i)} & \quad -({\rm id}_{K_1} \otimes \varphi_2)(\xi_1(\lambda)) +
        (\partial_2 \otimes {\rm id}_{F_1})(\xi_2(\lambda)) = 0\,, \\
\mbox{(ii)} & \quad ({\rm id}_{K_2} \otimes \varphi_1)(\xi_2(\lambda)) +
        (\partial_3 \otimes {\rm id}_{F_0})(\xi_3(\lambda)) = 0\,. 
\end{align*}
Here we write
\[
\xi_2(\lambda) = -\sum_{i = 1}^3 \check{e}_i \otimes u_{(i, \lambda)}
\hspace{2ex}\mbox{and}\hspace{2ex}
\xi_3(\lambda) = (e_1 \wedge e_2 \wedge e_3) \otimes d_\lambda\,,
\]
where $u_{(i, \lambda)} \in F_1$ and $d_\lambda \in F_0 = R$\,.
Then by (i) we get
\[
{\rm (iii)} \hspace{2ex} 
\left\{\begin{array}{l}
\vspace{0.3em}
\varphi_2(v_{(1, \lambda)}) = -x_2 \!\cdot\! u_{(3, \lambda)} - x_3 \!\cdot\! u_{(2, \lambda)} \\
\vspace{0.3em}
\varphi_2(v_{(2, \lambda)}) = x_1 \!\cdot\! u_{(3, \lambda)} - x_3 \!\cdot\! u_{(1, \lambda)} \\
\varphi_2(v_{(3, \lambda)}) = x_1 \!\cdot\! u_{(2, \lambda)} + x_2 \!\cdot\! u_{(1, \lambda)}
\end{array}\right.
\]
for any $\lambda \in \Lambda$\,.
Moreover, the equality (ii) implies
\[
{\rm (iv)} \hspace{2ex}
x_i \!\cdot\! d_\lambda = (-1)^{i-1} \!\cdot\! \varphi_1(u_{(i, \lambda)})
\]
for any $(i, \lambda) \in \tilde{\Lambda}$\,.
Now we define $\sigma_0 : K_0 \otimes_R F_3 \ra F_0$ and
$\sigma_1 : K_1 \otimes_R F_3 \ra F_1$ by setting
$\sigma_0(1 \otimes w_\lambda) = d_\lambda$ for any $\lambda \in \Lambda$
and $\sigma_1(e_i \otimes w_\lambda) = (-1)^{i-1} \!\cdot\! u_{(i, \lambda)}$
for any $(i, \lambda) \in \tilde{\Lambda}$\,.
The maps $\sigma_2 : K_2 \otimes_R F_3 \ra F_2$ and
$\sigma_3 : K_3 \otimes_R F_3 \ra F_3$ are defined as in (3) and (4).
Then, by (iii) and (iv), we see that
$\sigma_\bullet : K_\bullet \otimes_R F_3 \ra F_\bullet$
is a chain map.

Let us prove (1).
Because ${\rm Im}\,(\partial_1 \otimes {\rm id}_{F_3})$
is obviously contained in $\sigma_0^{-1}({\rm Im}\,\varphi_1)$\,,
we have to show the converse inclusion.
Take any $\eta_0 \in K_0 \otimes_R F_3$ such that
$\sigma_0(\eta_0) \in {\rm Im}\,\varphi_1$\,.
We write
\[
\eta_0 = \sum_{\lambda \in \Lambda} a_\lambda \otimes w_\lambda = 
\sum_{\lambda \in \Lambda} a_\lambda \!\cdot\! \xi_0(\lambda)\,,
\]
where $a_\lambda \in R$ for any $\lambda \in \Lambda$\,,
and set
\[
\eta_i = \sum_{\lambda \in \Lambda} a_\lambda \!\cdot\! \xi_i(\lambda)
\in K_i \otimes_R F_{3-i}
\]
for $i = 1, 2, 3$\,.
Then,
\[
\trmatrix{\,(\,\eta_0 \hspace{1ex} \eta_1 \hspace{1ex}
                           \eta_2 \hspace{1ex} \eta_3\,)} =
\sum_{\lambda \in \Lambda}
a_\lambda \!\cdot\!
\trmatrix{\,(\,\xi_0(\lambda) \hspace{1ex} \xi_1(\lambda) \hspace{1ex}
                    \xi_2(\lambda) \hspace{1ex} \xi_3(\lambda)\,)}
\in {\rm Ker}\,d_3\,.
\]
Furthermore, we have
\begin{eqnarray*}
\eta_3 & = & \sum_{\lambda \in \Lambda}
             a_\lambda \!\cdot\! ((e_1 \wedge e_2 \wedge e_3) \otimes d_\lambda) \\
  & = & (e_1 \wedge e_2 \wedge e_3) \otimes 
        \sum_{\lambda \in \Lambda} a_\lambda \!\cdot\! \sigma_0(1 \otimes w_\lambda) \\
  & = & (e_1 \wedge e_2 \wedge e_3) \otimes \sigma_0(\eta_0)\,,
\end{eqnarray*}
and so $\eta_3 \in {\rm Im}\,({\rm id}_{K_3} \otimes \varphi_1)$\,.
Hence $\eta_0 \in {\rm Im}\,(\partial_1 \otimes {\rm id}_{F_3})$
by (1) of \ref{2.2}.

Finally we prove (2).
Let us consider the following commutative diagram
\[
\begin{array}{cccccccl}
K_1 \otimes_R F_3 & \homom{\partial_1 \otimes {\rm id}_{F_3}} &
   K_0 \otimes_R F_3 & \ra & F_3 / QF_3 & \ra & 0 & {\rm (ex)} \\
\downarrow\vcenter{\rlap{$\scriptstyle{\sigma_1}$}} &   &
  \downarrow\vcenter{\rlap{$\scriptstyle{\sigma_0}$}} &   & 
    \downarrow\vcenter{\rlap{$\scriptstyle{\overline{\sigma}_0}$}} &
      &   &   \\
F_1 & \homom{\varphi_1} & F_0 & \ra & R / \ga & \ra & 0 & {\rm (ex)\,.}
\end{array}
\]
where $\overline{\sigma}_0$ is the map indeuced from $\sigma_0$\,.
Notice that by (iv) we have $d_\lambda \in \ga :_R Q$
for any $\lambda \in \Lambda$\,,
and so ${\rm Im}\,\sigma_0 \subseteq \ga :_R Q$\,.
Hence ${\rm Im}\,\overline{\sigma}_0 \subseteq (\ga :_R Q) / \ga$\,.
On the other hand, as $\sigma_0^{-1}({\rm Im}\,\varphi_1) =
{\rm Im}\,(\partial_1 \otimes {\rm id}_{F_3})$\,,
we see that $\overline{\sigma}_0$ is injective.
Therefore we get ${\rm Im}\,\overline{\sigma}_0 = (\ga :_R Q) / \ga$
since $(\ga :_R Q) / \ga \cong F_3 / QF_3$ by \ref{3.1}
and $F_3 / QF_3$ has a finite length.
Thus the assertion (2) follows and the proof is completed.

\vspace{1em}

In the rest,
$\sigma_\bullet : K_\bullet \otimes_R F_3 \ra F_\bullet$
is the chain map constructed in \ref{3.2}.
Then, by \ref{2.1} the mapping cone ${\rm Con}\,(\sigma_\bullet)$
gives an $R$-free resolution of $\ga :_R Q$\,,
that is,
\[
0
\shomom{}
K_3 \otimes_R F_3
\homom{\psi_4}
\begin{array}{c}
K_2 \otimes_R F_3 \\
\oplus \\
F_3
\end{array}
\homom{\psi_3}
\begin{array}{c}
K_1 \otimes_R F_3 \\
\oplus \\
F_2
\end{array}
\homom{\lsp{\varphi_2}}
\begin{array}{c}
K_0 \otimes_R F_3 \\
\oplus \\
F_1
\end{array}
\homom{\lss{\varphi_1}}
F_0 = R\,, 
\]
is acyclic and ${\rm Im}\,\lss{\varphi_1} = \ga :_R Q$\,, where
\[
\psi_4 = 
\left(\begin{array}{c}
\partial_3 \otimes {\rm id}_{F_3} \\
-\sigma_3
\end{array}\right),\,
\psi_3 =
\left(\begin{array}{cc}
\partial_2 \otimes {\rm id}_{F_3} & 0 \\
\sigma_2 & \varphi_3
\end{array}\right),\,
\lsp{\varphi_2} =
\left(\begin{array}{cc}
\partial_1 \otimes {\rm id}_{F_3} & 0 \\
-\sigma_1 & \varphi_2
\end{array}\right),\,
\lss{\varphi_1} =
( \sigma_0 \hspace{1ex} \varphi_1 ) \,.
\]
Because $\sigma_3 : K_3 \otimes_R F_3 \ra F_3$ is an isomorphism,
we can define
\[
\phi = (\,0 \hspace{1.5ex} -\sigma_3^{-1}\,) \hspace{1ex} : 
\begin{array}{c}
K_2 \otimes_R F_3 \\
\oplus \\
F_3
\end{array}
\ra \hspace{1ex} K_3 \otimes_R F_3\,.
\]
Then $\phi \circ \psi_4 = {\rm id}_{K_3 \otimes_R F_3}$ and
${\rm Ker}\,\phi = K_2 \otimes_R F_3$\,.
Hence, by (1) of \ref{2.3}, we get the acyclic complex
\[
0 \ra \lsp{F_3} \homom{\lsp{\varphi_3}} \lsp{F_2} \homom{\lsp{\varphi_2}} 
\lss{F_1} \homom{\lss{\varphi_1}} \lss{F_0} = R\,, 
\]
where
\[
\lsp{F_3} = K_2 \otimes_R F_3\,, \hspace{2ex}
\lsp{F_2} = 
\begin{array}{c}
K_1 \otimes_R F_3 \\
\oplus \\
F_2
\end{array}\,, \hspace{2ex}
\lss{F_1} = 
\begin{array}{c}
K_0 \otimes_R F_3 \\
\oplus \\
F_1
\end{array}\hspace{1ex} \mbox{and} \hspace{2ex}
\lsp{\varphi_3} = \left(
\begin{array}{c}
\partial_2 \otimes {\rm id}_{F_3} \\
\sigma_2
\end{array}\right)\,. 
\]

Although ${\rm Im}\,\lsp{\varphi_3}$ may not be contained in $\gm\,\lsp{F_2}$\,,
removing non-minimal components from $\lsp{F_3}$ and $\lsp{F_2}$\,,
we get free $R$-modules $\lss{F_3}$ and $\lss{F_2}$ such that
\[
0 \ra \lss{F_3} \homom{\lss{\varphi_3}} \lss{F_2}
\homom{\lss{\varphi_2}} \lss{F_1} \homom{\lss{\varphi_1}}
\lss{F_0} = R
\]
is acyclic and ${\rm Im}\,\lss{\varphi_3} \subseteq \gm\,\lss{F_2}$\,,
where $\lss{\varphi_3}$ and $\lss{\varphi_2}$ are the restrictions of
$\lsp{\varphi_3}$ and $\lsp{\varphi_2}$\,, respectively.
In the rest of this section,
we describe a concrete procedure to get $\lss{F_3}$ and $\lss{F_2}$\,.
For that purpose, we use the following notation.
As is described in Introduction,
for any $\xi \in K_1 \otimes_R F_3$ and $\eta \in F_2$\,,
\[
\lrbr{\xi} := \binom{\xi}{0} \in \lsp{F_2} 
\hspace{3ex} \mbox{and} \hspace{3ex}
\lran{\eta} := \binom{0}{\eta} \in \lsp{F_2}\,.
\]
In particular, for any $(i, \lambda) \in \tilde{\Lambda}$\,, we denote
$\lrbr{e_i \otimes w_\lambda}$ by $[i , \lambda]$\,.
Moreover, for a subset $U \subseteq F_2$\,,
$\lran{U} := \{ \lran{u} \}_{u \in U}$\,.

Now, let us choose a subset $\lsp{\Lambda}$ of $\tilde{\Lambda}$
and a subset $U$ of $F_2$ so that
\[
\{ v_{(i , \lambda)} \}_{(i , \lambda)\, \in\, \lsp{\Lambda}} \cup U
\]
is an $R$-free basis of $F_2$\,.
We would like to choose $\lsp{\Lambda}$ as big as possible.
The following almost obvious fact is useful to find $\lsp{\Lambda}$ and $U$\,. 

\begin{lem}\label{3.3}
Let $V$ be an $R$-free basis of $F_2$\,.
If a subset $\lsp{\Lambda}$ of $\tilde{\Lambda}$ and a subset $U$ of $V$ satisfy
\begin{itemize}
\item[{\rm (i)}]
$|\,\lsp{\Lambda}\,| + |\,U\,| \leq |\,V\,|$\,,
where $|\,\cdot\,|$ denotes the number of  elements of the set, and
\item[{\rm (ii)}]
$V \subseteq R \!\cdot\! \{\,v_{(i, \lambda)}\,\}_{(i, \lambda)\,\in\,\lsp{\Lambda}} +
R \!\cdot\! U + \gm F_2$\,,
\end{itemize}
then $\{\,v_{(i , \lambda)}\,\}_{(i , \lambda)\, \in\, \lsp{\Lambda}} \cup U$
is an $R$-free basis of $F_2$\,.
\end{lem}

\noindent
Let us notice that
\[
\{ [i, \lambda] \}_{(i, \lambda)\,\in\,\tilde{\Lambda}} \,\cup\,
\{ \lran{v_{(i, \lambda)}} \}_{(i, \lambda)\,\in\,\lsp{\Lambda}} \cup \lran{U}
\]
is an $R$-free basis of $\lsp{F_2}$\,.
We define $\lss{F_2}$ to be the direct summand of $\lsp{F_2}$
generated by
\[
\{ [i, \lambda] \}_{(i, \lambda)\,\in\,\tilde{\Lambda}}
\,\cup\, \lran{U}\,.
\]
Let $\lss{\varphi_2} $ be the restriction of $\lsp{\varphi_2}$ to $\lss{F_2}$\,.
We need the next result at the final step of the process to
compute symbolic powers.

\begin{thm}\label{3.4}
If we can take $\tilde{\Lambda}$ itself as $\lsp{\Lambda}$\,, then
\[
0 \ra \lss{F_2} \homom{\lss{\varphi_2}} \lss{F_1}
\homom{\lss{\varphi_1}} \lss{F_0} = R
\]
is acyclic. Hence we have
$\dep{}{R / (\ga :_R Q)} > 0$\,.
\end{thm}

\noindent
{\it Proof}.\,
If $\lsp{\Lambda} = \tilde{\Lambda}$\,,
there exists a homomorphism 
$\phi : \lsp{F_2} \ra \lsp{F_3}$ such that
\begin{eqnarray*}
\phi\,([i, \lambda]) & = &
  0 \hspace{2ex} \mbox{for any $(i, \lambda) \in \tilde{\Lambda}$\,,}  \\
\phi\,(\lran{v_{(i, \lambda)}}) & = &
  (-1)^i\!\cdot\!\check{e}_i \otimes w_\lambda
  \hspace{2ex} \mbox{for any $(i, \lambda) \in \tilde{\Lambda}$\,,}  \\
\phi\,(\lran{u}) & = & 
  0 \hspace{2ex} \mbox{for any $u \in U$\,.}
\end{eqnarray*}
Then $\phi\circ\lsp{\varphi_3} = {\rm id}_{\lsp{F_3}}$
and ${\rm Ker}\,\phi = \lss{F_2}$\,.
Hence by (1) of \ref{2.3} we get the required assertion.

\vspace{2em}
In the rest of this section, we assume
$\lsp{\Lambda} \subsetneq \tilde{\Lambda}$
and put $\lss{\Lambda} = \tilde{\Lambda} \setminus \lsp{\Lambda}$\,.
Then, for any $(j, \mu) \in \lss{\Lambda}$\,,
it is possible to write
\[
v_{(j, \mu)} = \sum_{(i, \lambda)\,\in\,\lsp{\Lambda}}
a_{(i, \lambda)}^{(j, \mu)}\!\cdot\! v_{(i, \lambda)} +
\sum_{u \in U} b_u^{(j, \mu)}\!\cdot\! u\,,
\]
where $a_{(i, \lambda)}^{(j, \mu)}\,, b_u^{(j, \mu)} \in R$\,.
Here, if $\lsp{\Lambda}$ is big enough,
we can choose every $b_u^{(j, \mu)}$ from $\gm$\,.
In fact, if $b_u^{(j, \mu)} \not\in \gm$ for some $u \in U$\,,
then we can replace $\lsp{\Lambda}$ and $U$ by
$\lsp{\Lambda} \cup \{\,(j, \mu)\,\}$ and 
$U \setminus \{\,u\,\}$\,, respectively.
Furthermore, because of a practical reason,
let us allow that some terms of $v_{(i, \lambda)}$
for $(i, \lambda) \in \lss{\Lambda}$ with non-unit coefficients
appear in the right hand side, that is, 
for any $(j, \mu) \in \lss{\Lambda}$\,, we write
\[
v_{(j, \mu)} = \sum_{(i, \lambda)\,\in\,\tilde{\Lambda}}
a_{(i, \lambda)}^{(j, \mu)}\!\cdot\! v_{(i, \lambda)} +
\sum_{u \in U} b_u^{(j, \mu)}\!\cdot\! u\,,
\]
where
\[
a_{(i, \lambda)}^{(j, \mu)} \in \left\{
\begin{array}{ll}
R & \mbox{if $(i, \lambda) \in \lsp{\Lambda}$} \\
\gm & \mbox{if $(i, \lambda) \in \tilde{\Lambda}$}
\end{array}\right.
\hspace{2ex}\mbox{and}\hspace{2ex}
b_u^{(j, \mu)} \in \gm\,.
\]
Using this expression, for any $(j, \mu) \in \lss{\Lambda}$\,,
the following element in $\lsp{F_3}$ can be defined.
\[
\lss{w_{(j, \mu)}} :=
(-1)^j\!\cdot\!\check{e}_j \otimes w_\mu + 
\sum_{(i, \lambda)\,\in\,\tilde{\Lambda}}
(-1)^{i+1}a_{(i, \lambda)}^{(j, \mu)}\!\cdot\!
\check{e}_i \otimes w_\lambda\,.
\]

\begin{lem}\label{3.5}
For any $i \in \{ 1, 2, 3 \}$\,,
let $s_i$ and $t_i$ be integers such that $s_i < t_i$ and
$\{ 1, 2, 3 \} = \{ i, s_i, t_i \}$\,.
Then, for any $(j, \mu) \in \lss{\Lambda}$\,, we have
\begin{equation*}
\begin{split}
\lsp{\varphi_3}(\lss{w_{(j, \mu)}})
 &= (-1)^jx_{s_j}\!\cdot\![t_j, \mu] + 
    (-1)^{j+1}x_{t_j}\!\cdot\![s_j, \mu] \\
 & \quad + \sum_{(i, \lambda)\,\in\,\tilde{\Lambda}}
        (-1)^{i+1}a_{(i, \lambda)}^{(j, \mu)}\!\cdot\!
        \{ x_{s_i}\!\cdot\![t_i, \lambda] - x_{t_i}\!\cdot\![s_i, \lambda] \}
      + \sum_{u \in U} b_u^{(j, \mu)}\!\cdot\!\lran{u}\,.
\end{split}
\end{equation*}
As a consequence, we have
$\lsp{\varphi_3}(\lss{w_{(j, \mu)}}) \in \gm\,\lss{F_2}$
for any $(j, \mu) \in \lss{\Lambda}$\,.
\end{lem}

\noindent
{\it Proof}.\,
By the definition of $\lsp{\varphi_3}$\,,
for any $(j, \mu) \in \lss{\Lambda}$\,, we have
\[
\lsp{\varphi_3}(\lss{w_{(j, \mu)}}) = 
\lrbr{(\partial_2 \otimes {\rm id}_{F_3})(\lss{w_{(j, \mu)}})} +
\lran{\sigma_2(\lss{w_{(j, \mu)}})}\,.
\]
Because
\begin{equation*}
\begin{split}
(\partial_2 \otimes {\rm id}_{F_3})(\lss{w_{(j, \mu)}})
 &=
(-1)^j\!\cdot\!\partial_2\check{e}_j \otimes w_\mu +
\sum_{(i, \lambda)\,\in\,\tilde{\Lambda}}
(-1)^{i+1}a_{(i, \lambda)}^{(j, \mu)}\!\cdot\!
\partial_2\check{e}_i \otimes w_\lambda  \\
  &=
(-1)^j\!\cdot\!\partial_2(e_{s_j}\wedge e_{t_j}) \otimes w_\mu +
\sum_{(i, \lambda)\,\in\,\tilde{\Lambda}}
(-1)^{i+1}a_{(i, \lambda)}^{(j, \mu)}\!\cdot\!
\partial_2(e_{s_i}\wedge e_{t_i}) \otimes w_\lambda \\
  &=
(-1)^j\!\cdot\! (x_{s_j}e_{t_j} - x_{t_j}e_{s_j}) \otimes w_\mu  \\
  & \hspace*{15ex}
+ \sum_{(i, \lambda)\,\in\,\tilde{\Lambda}}
(-1)^{i+1}a_{(i, \lambda)}^{(j, \mu)}\!\cdot\!
(x_{s_i}e_{t_i} - x_{t_i}e_{s_i}) \otimes w_\lambda
\end{split}
\end{equation*}
and
\begin{eqnarray*}
\sigma_2(\lss{w_{(j, \mu)}}) & = &
(-1)^j\!\cdot\!\sigma_2(\check{e}_j \otimes w_\mu) +
\sum_{(i, \lambda)\,\in\,\tilde{\Lambda}}
(-1)^{i+1}a_{(i, \lambda)}^{(j, \mu)}\!\cdot\!
\sigma_2(\check{e}_i \otimes w_\lambda) \\
 & = &
v_{(j, \mu)} - \sum_{(i, \lambda)\,\in\,\tilde{\Lambda}}
a_{(i, \lambda)}^{(j, \mu)}\!\cdot\! v_{(i, \lambda)} \\
 & = &
\sum_{u \in U} b_u^{(j, \mu)}\!\cdot\! u\,,
\end{eqnarray*}
we get the required equality.

\vspace{2ex}
Let $\lss{F_3}$ be the $R$-submodule of $\lsp{F_3}$
generated by $\{\,\lss{w_{(j, \mu)}}\,\}_{(j, \mu)\,\in\,\lss{\Lambda}}$
and let $\lss{\varphi_3}$ be the restriction of $\lsp{\varphi_3}$
to $\lss{F_3}$\,.
By \ref{3.5} we have ${\rm Im}\,\lss{\varphi_3} \subseteq \lss{F_2}$\,.
Thus we get a complex
\[
0 \ra \lss{F_3} \homom{\lss{\varphi_3}} \lss{F_2}
\homom{\lss{\varphi_2}} \lss{F_1} \homom{\lss{\varphi_1}}
\lss{F_0} = R\,.
\]
This is the complex we desire.
In fact, the following result holds.

\begin{thm}\label{3.6}
$(\lss{F_\bullet}\,,\,\lss{\varphi_\bullet})$ is an acyclic complex of finitely generated
free $R$-modules with the following properties.
\begin{itemize}
\item[{\rm (1)}]
${\rm Im}\,\lss{\varphi_1} = \ga :_R Q$\,and
${\rm Im}\,\lss{\varphi_3} \subseteq \gm\,\lss{F_2}$\,.
\item[{\rm (2)}]
$\{\,\lss{w_{(j, \mu)}}\,\}_{(j, \mu)\,\in\,\lss{\Lambda}}$
is an $R$-free basis of $\lss{F_3}$\,.
\item[{\rm (3)}]
$\{\,[i, \lambda]\,\}_{(i, \lambda)\,\in\,\tilde{\Lambda}}
\,\cup\, \lran{U}$
is an $R$-free basis of $\lss{F_2}$\,.
\end{itemize}
\end{thm}

\noindent
{\it Proof}.\,
First, let us notice that
$\{ \check{e}_i \otimes w_\lambda \}_{(i, \lambda)\,\in\,\tilde{\Lambda}}$
is an $R$-free basis of $\lsp{F_3}$ and
\[
\check{e}_j \otimes w_\mu \in R\!\cdot\!\lss{w_{(j, \mu)}} +
R\!\cdot\!
\{ \check{e}_i \otimes w_\lambda \}_{(i, \lambda)\,\in\,\lsp{\Lambda}} +
\gm\,\lsp{F_3}
\]
for any $(j, \mu) \in \lss{\Lambda}$\,.
Hence, by Nakayama's lemma it follows that $\lsp{F_3}$ is generated by
\[
\{ \check{e}_i \otimes w_\lambda \}_{(i, \lambda)\,\in\,\lsp{\Lambda}}
\,\cup\,
\{ \lss{w_{(j, \mu)}} \}_{(j, \mu)\,\in\,\lss{\Lambda}}\,,
\]
which must be an $R$-free basis since
${\rm rank}_R\,\lsp{F_3} = |\,\tilde{\Lambda}\,| =
|\,\lsp{\Lambda}\,| + |\,\lss{\Lambda}\,|$\,.
Let $\lspp{F_3}$ be the $R$-submodule of $\lsp{F_3}$ generated by
$\{\,\check{e}_i \otimes w_\lambda\,\}_{(i, \lambda)\,\in\,\lsp{\Lambda}}$\,.
Then $\lspp{F_3} \oplus \lss{F_3}$\,.

Next, let us recall that
\[
\{ [i, \lambda] \}_{(i, \lambda)\,\in\,\tilde{\Lambda}}
\,\cup\,
\{ \lran{v_{(i, \lambda)}} \}_{(i, \lambda)\,\in\,\lsp{\Lambda}}
\cup \lran{U}
\]
is an $R$-free basis of $\lsp{F_2}$\,.
Because
\[
\lran{v_{(i, \lambda)}} = 
(-1)^i\!\cdot\!\lsp{\varphi_3}(\check{e}_i \otimes w_\lambda) +
(-1)^i\!\cdot\!\lrbr{\partial_2\check{e}_i \otimes w_\lambda}\,,
\]
we see that
\[
\{ [i, \lambda] \}_{(i, \lambda)\,\in\,\tilde{\Lambda}}
\,\cup\,
\{ \lsp{\varphi_3}
(\check{e}_i \otimes w_\lambda)  \}_{(i, \lambda)\,\in\,\lsp{\Lambda}}
\,\cup\, \lran{U}
\]
is also an $R$-free basis.
Let $\lspp{F_2} = R\!\cdot\!
\{\,\lsp{\varphi_3}
(\check{e}_i \otimes w_\lambda) \,\}_{(i, \lambda)\,\in\,\lsp{\Lambda}}$\,.
Then $\lsp{F_2} = \lspp{F_2} \oplus \lss{F_2}$\,.

It is obvious that
$\lsp{\varphi_3}(\lspp{F_3}) = \lspp{F_2}$\,.
Moreover, by \ref{3.5} we get
$\lsp{\varphi_3}(\lss{F_3}) \subseteq \lss{F_2}$\,.
Therefore, by (2) of \ref{2.3}, it follows that
$\lss{F_\bullet}$ is acyclic.
We have already seen (3) and the first assertion of (1).
The second assertion of (1) follows from \ref{3.5}.
Moreover, the assertion (2) is now obvious.

\section{Computing symbolic powers}

Let $x, y, z$ be an sop for $R$ and $I$ an ideal of $R$
generated by the maximal minors of the matrix
\[
\Phi = \left(\begin{array}{lll}
x^\alpha & y^\beta & z^\gamma \\
y^{\beta'} & z^{\gamma'} & x^{\alpha'}
\end{array}\right)\,,
\]
where $\alpha, \beta, \gamma, \alpha', \beta', \gamma'$
are positive integers.
As is well known,
$R / I$ is a Cohen-Macaulay ring with $\dim R / I = 1$\,.
In this section,
we give a minimal free resolution of $I^n$ for any $n > 0$\,,
and consider its $\ast$-transform in order to
compute the symbolic power $I^{(n)}$\,.
We put
\[
a = z^{\gamma + \gamma'} - x^{\alpha'}y^\beta\,,\,
b = x^{\alpha + \alpha'} - y^{\beta'}z^\gamma\,,\,
c = y^{\beta + \beta'} - x^\alpha z^{\gamma'}\,.
\]
Then, $I = (a, b, c)R$ and we have the next result
(See \cite{h1} for the definition of d-sequences).

\begin{lem}\label{4.1}
The following assertions hold.
\begin{itemize}
\item[{\rm (1)}]
$x^\alpha a + y^\beta b + z^\gamma c = 0$ and
$y^{\beta'}a + z^{\gamma'}b + x^{\alpha'}c = 0$\,.
\item[{\rm (2)}]
Let $\gp \in \ass{R}{R / I}$\,.
Then $IR_\gp$ is generated by any two elements of
$a, b, c$\,.
\item[{\rm (3)}]
Any two elements of $a, b, c$ form an ssop for $R$\,.
\item[{\rm (4)}]
$a, b, c$ is an unconditioned d-sequence.
\end{itemize}
\end{lem}

\noindent
{\it Proof}.\,
(1)\,
These equalities can be checked directly.

(2)\,
Let us prove $IR_\gp = (a, b)R_\gp$\,.
If $x \in \gp$\,,
then $y, z \in \sqrt{(a, c, x)R} \subseteq \gp$\,,
and so $\gp = \gm$\,,
which contradicts to the Cohen-Macaulayness of $R / I$\,.
Hence $x \not\in \gp$\,.
Then 
\[
c = -(y^{\beta'} a + z^{\gamma'} b) / x^{\alpha'} \in (a, b)R_\gp\,,
\]
which means $IR_\gp = (a, b)R_\gp$\,.

(3)\,
For example,
as $x, z \in \sqrt{(a, b, y)R}$\,,
it follows that $a, b$ is an ssop for $R$\,.

(4)\,
Let us prove that $a, b, c$ is a d-sequence.
As $a, b$ is a regular sequence by (3),
it is enough to show $(a, b)R :_R c^2 = (a, b)R :_R c$\,.
We obviously have $(a, b)R :_R c^2 \supseteq (a, b)R :_R c$\,.
Take any $\gq \in \ass{R}{R / (a, b)R :_R c}$\,.
As $R / (a, b)R :_R c \hookrightarrow R / (a, b)R$\,,
we have $\height{R}{\gq} = 2$\,.
If $c \in \gq$\,,
then $\gq \in \Min{R}{R / I}$\,,
and so $IR_\gq = (a, b)R_\gq$ by (2),
which means
\[
(a, b)R_\gq :_{R_\gq} c^2 = (a, b)R_\gq :_{R_\gq} c = R_\gq\,.
\]
If $c \not\in \gq$\,, we have
\[
(a, b)R_\gq :_{R_\gq} c^2 = (a, b)R_\gq :_{R_\gq} c = (a, b)R_\gq\,.
\]
Therefore we get the required equality.

\vspace{1em}
We take an indeterminate $t$ and consider the Rees algebra $R[ It ]$\,.
Moreover, we take three indeterminates $A, B, C$
and put $S = R[A, B, C]$\,.
We regard $S$ as a $\zz$-graded ring by setting
$\deg A = \deg B = \deg C = 1$\,.
Let $\pi : S \ra R[ It ]$ be the graded homomorphism
of $R$-algebras such that $\pi(A) = at$, $\pi(B) = bt$ and $\pi(C) = ct$\,.
By (4) of \ref{4.1} it follows that ${\rm Ker}\,\pi$ is generated by
linear forms (cf. \cite[Theorem 3.1]{h2}).
On the other hand,
\[
0 \ra R^{\oplus 2} \homom{\trmatrix{\,\Phi}} R^{\oplus 3}
\homom{(a\hspace{1ex}b\hspace{1ex}c)} R \ra R / I \ra 0
\]
is a minimal free resolution of $R / I$\,.
Hence ${\rm Ker}\,\pi$ is generated by
\[
f := x^\alpha A + y^\beta B + z^\gamma C 
\hspace{2ex}\mbox{and}\hspace{2ex}
g := y^{\beta'}A + z^{\gamma'}B + x^{\alpha'}C\,.
\]
Thus we get $S / (f, g)S \cong R[ It ]$\,.
Then, as $f, g$ is a regular sequence of $S$\,,
\[
0 \ra S(-2) \homom{\binom{-g}{f}} S(-1) \oplus S(-1)
\homom{(f\hspace{1ex}g)} S \homom{\pi} R[ It ] \ra 0
\]
is a graded $S$-free resolution of $R[ It ]$\,.
Now we take its homogeneous part of degree $n$\,,
and get the next result.

\begin{thm}\label{4.2}
For any $n \geq 2$\,,
\[
0 \ra S_{n-2} \homom{\binom{-g}{f}} S_{n-1} \oplus S_{n-1}
\homom{(f\hspace{1ex}g)} S_n \homom{\epsilon} R
\]
is acyclic and it is a minimal free resolution of $I^n$\,,
where $S_d$ {\rm ($d \in \zz$)} is the $R$-submodule of $S$
consisting of homogeneous elements of degree $d$ and
$\epsilon$ is the $R$-linear map defined by substituting
$a, b, c$ for $A, B, C$\,, respectively.
\end{thm}

Let us denote the complex in \ref{4.2} by $(F_\bullet^1, \varphi_\bullet^1)$\,,
that is, we set
\[
\begin{array}{l}
\vspace{0.3em}
F_3^1 = S_{n-2}\,,\,F_2^1 = S_{n-1} \oplus S_{n-1}\,,\,F_1^1 = S_n\,,\,F_0^1 = R\,, \\
\varphi_3^1 = \binom{-g}{f}\,,\,\varphi_2^1 = (f\hspace{1ex}g)
\hspace{2ex}\mbox{and}\hspace{2ex}
\varphi_1^1 = \epsilon\,.
\end{array}
\]
Then $F_\bullet^1$ is an acyclic complex of finitely generated free $R$-modules
and ${\rm Im}\,\varphi_1^1 = I^n$\,.
The number "$1$" of  $F_\bullet^1$ means that
it is the first acyclic complex we need for computing $I^{(n)}$\,.
Our strategy is as follows.
Taking the $\ast$-transform of $F^1_\bullet$ with respect to
suitable powers of $x, y, z$\,, we get 
$\lss{F_\bullet^1}$\,,
which is denoted by
$F_\bullet^2$\,.
If its length is still $3$,
we again take some $\ast$-transform of $F_\bullet^2$
and get $F_\bullet^3$\,.
By repeating this operation successively,
we eventually get an acyclic complex
$F_\bullet^k$ of length $2$\,.
Then the family $\{\,F^i_\bullet\,\}_{1 \leq i \leq k}$
of acyclic complexes has 
the complete information on $I^{(n)}$\,.

Let $\alpha'' = \min\{\alpha, \alpha'\}$\,,
$\beta'' = \min\{\beta\,, \beta'\}$\,,
$\gamma'' = \min\{\gamma\,, \gamma'\}$ and
$Q = (x^{\alpha''},\, y^{\beta''},\, z^{\gamma''})R$\,.
Because $f$ and $g$ are elements of $Q$\,,
we have ${\rm Im}\,\varphi_3^1 \subseteq QF_2^1$\,,
and so by \ref{3.1} we get the following.

\begin{thm}\label{4.3}
$(I^n :_R Q) / I^n \cong (R / Q)^{\oplus \binom{n}{2}}$\,.
\end{thm}

Now we are going to take the $\ast$-transform of $F_\bullet^1$
with respect to $x^{\alpha''},\, y^{\beta''},\, z^{\gamma''}$\,.
At first, we have to fix $\Lambda^{\!1}$\,,
which is an $R$-free basis of $F^1_3$\,.
For any $0 \leq d \in \zz$\,,
let us denote by $\mon{A, B, C}{d}$ the set
$\{ A^iB^jC^k \mid 0 \leq i, j, k \in \zz
\hspace{1ex}\mbox{and}\hspace{1ex}
i + j + k = d \}$\,,
which is an $R$-free basis of $S_d$\,.
We take $\mon{A, B, C}{n - 2}$ as $\Lambda^{\!1}$\,.
Then, for any $M \in \mon{A, B, C}{n - 2}$\,,
we have to write
\[
\varphi^1_3(M) =
x^{\alpha''}\!\cdot\! v^1_{(1, M)} +
y^{\beta''}\!\cdot\! v^1_{(2, M)} +
z^{\gamma''}\!\cdot\! v^1_{(3, M)}\,,
\]
where $v^1_{(i, M)} \in F^1_2$ for $i = 1, 2, 3$\,.
As is described at the end of Introduction,
for $h \in S_{n-1}$\,,
let us denote the elements
\[
\binom{h}{0}
\hspace{1ex}\mbox{and}\hspace{1ex}
\binom{0}{h}
\in F^1_2 = S_{n-1} \oplus S_{n-1}
\]
by $\lrbr{h}$ and $\lran{h}$\,, respectively.
Then, for any $M \in \mon{A, B, C}{n - 2}$\,, we have
\begin{equation*}
\begin{split}
\varphi^1_3(M) &= \binom{-gM}{fM} \\
 &= 
\binom{-y^{\beta'}AM - z^{\gamma'}BM -
x^{\alpha'}CM}{x^\alpha AM + y^\beta BM + z^\gamma CM} \\
 &=
-y^{\beta'}\!\cdot\!\lrbr{AM} - z^{\gamma'}\!\cdot\!\lrbr{BM} -
x^{\alpha'}\!\cdot\!\lrbr{CM}  \\
 & \hspace{10ex}
+ x^\alpha\!\cdot\!\lran{AM} +
y^\beta\!\cdot\!\lran{BM} + z^\gamma\!\cdot\!\lran{CM} \\
 &=
x^{\alpha''}\!\cdot\! v^1_{(1, M)} +
y^{\beta''}\!\cdot\! v^1_{(2, M)} +
z^{\gamma''}\!\cdot\! v^1_{(3, M)}\,,
\end{split}
\end{equation*}
where
\begin{eqnarray*}
v^1_{(1, M)} & = &
x^{\alpha - \alpha''}\!\cdot\!\lran{AM} -
x^{\alpha' - \alpha''}\!\cdot\!\lrbr{CM}\,,\, \\
v^1_{(2, M)} & = &
y^{\beta - \beta''}\!\cdot\!\lran{BM} -
y^{\beta' - \beta''}\!\cdot\!\lrbr{AM}\,,\, \\
v^1_{(3, M)} & = &
z^{\gamma - \gamma''}\!\cdot\!\lran{CM} -
z^{\gamma' - \gamma''}\!\cdot\!\lrbr{BM}\,. \\
\end{eqnarray*}

We set $\widetilde{\Lambda^{\!1}} = \{1, 2, 3\} \times \Lambda^{\!1}$
and we have to choose its subset
$\lsp{\Lambda^{\!1}}$ as big as possible so that
$\{v^1_{(i, M)}\}_{(i, M)\,\in\,\lsp{\Lambda^{\!1}}}$
is a part of an $R$-free basis of $F^1_2$\,.
For that purpose, we need to fix a canonical $R$-free basis of $F^1_2$\,.
For a subset $H$ of $S_{n-1}$\,,
we denote the families
$\{ \lrbr{h} \}_{h \in H}$ and
$\{ \lran{h} \}_{h \in H}$ by
$\lrbr{H}$ and $\lran{H}$\,, respectively.
Let us notice that
$\lrbr{\mon{A, B, C}{n - 1}} \cup \lran{\mon{A, B, C}{n - 1}}$
is an $R$-free basis of $F^1_2$\,.

Setting $n = 2$\,, we get the next result.

\begin{thm}\label{4.4}
{\rm (cf. \cite{gns})}\,
$I^{(2)} = I^2 :_R Q$ and 
$I^{(2)} / I^2 \cong R / Q$\,.
\end{thm}

\noindent
{\it Proof}.\,
By replacing rows and columns of $\Phi$ and by replacing $x, y, z$
if necessary,
we may assume that one of the following conditions are satisfied;
\[
{\rm (i)}\,
\alpha \leq \alpha'\,,\, \beta \leq \beta'\,,\, \gamma \leq \gamma'\, ;\hspace{2ex}
{\rm (ii)}\,
\alpha \geq \alpha'\,,\, \beta \leq \beta'\,,\, \gamma \leq \gamma'\,.
\]
Let $n = 2$\,.
Then $\Lambda^{\!1} = \{\,1\,\}$\,.
In the case (i), we have
$\alpha'' = \alpha, \beta'' = \beta, \gamma'' = \gamma$ and
\[
v^1_{(1, 1)} = \lran{A} - x^{\alpha' - \alpha}\!\cdot\!\lrbr{C}\,,\,
v^1_{(2, 1)} = \lran{B} - y^{\beta' - \beta}\!\cdot\!\lrbr{A}\,,\,
v^1_{(3, 1)} = \lran{C} - z^{\gamma' - \gamma}\!\cdot\!\lrbr{B}\,.
\]
Then, as $\lrbr{\mon{A, B, C}{1}} \cup \lran{\mon{A, B, C}{1}}$
is an $R$-free basis, by \ref{3.3} we see that
$\{ v^1_{(i, 1)} \}_{i = 1, 2, 3} \cup \lrbr{\mon{A, B, C}{1}}$
is an $R$-free basis of $F^1_2$\,.
On the other hand, in the case (ii), we have
$\alpha'' = \alpha', \beta'' = \beta, \gamma'' = \gamma$ and
\[
v^1_{(1, 1)} = x^{\alpha - \alpha'}\!\cdot\!\lran{A} - \lrbr{C}\,,\,
v^1_{(2, 1)} = \lran{B} - y^{\beta' - \beta}\!\cdot\!\lrbr{A}\,,\,
v^1_{(3, 1)} = \lran{C} - z^{\gamma' - \gamma}\!\cdot\!\lrbr{B}\,.
\]
Then, $\{ v^1_{(i, 1)} \}_{i = 1, 2, 3} \cup
\{ \lran{A}, \lrbr{A}, \lrbr{B} \}$ is an $R$-free basis of $F^1_2$\,.
In either case, we can take
$\widetilde{\Lambda^{\!1}}$ as $\lsp{\Lambda^{\!1}}$\,.
Hence, by \ref{3.4} we see $\dep{}{R / (I^2 :_R Q)} > 0$\,,
and so $I^{(2)} = I^2 :_R Q$\,.
The second assertion follows from \ref{4.3}.

\vspace{1em}
Similarly as the proof of \ref{4.4},
in order to study $I^{(n)}$ for $n \geq 3$\,,
we have to consider dividing the situation into several cases.
In the rest of this section, let us assume
\[
\alpha = 1\,,\, \alpha' = 2\,,\, 2\beta \leq \beta'\,,\,
2\gamma \leq \gamma'\,,
\]
and explain how to compute $I^{(3)}$
using $\ast$-transforms.
We have
$\alpha'' = 1, \beta'' = \beta$ and $\gamma'' = \gamma$\,.
Let $n = 3$\,.
Then $\Lambda^{\!1} = \{A, B, C\}$ and we have
\begin{alignat*}{3}
v^1_{(1, A)} & = \lran{A^2} - x\!\cdot\!\lrbr{AC}\,,
& \hspace{1ex}
v^1_{(2, A)} & = \lran{AB} - y^{\beta' - \beta}\!\cdot\!\lrbr{A^2}\,,
& \hspace{1ex}
v^1_{(3, A)} & = \lran{AC} - z^{\gamma' - \gamma}\!\cdot\!\lrbr{AB}\,,
\\
v^1_{(1, B)} & = \lran{AB} - x\!\cdot\!\lrbr{BC}\,,
& \hspace{1ex}
v^1_{(2, B)} & = \lran{B^2} - y^{\beta' - \beta}\!\cdot\!\lrbr{AB}\,,
& \hspace{1ex}
v^1_{(3, B)} & = \lran{BC} - z^{\gamma' - \gamma}\!\cdot\!\lrbr{B^2}\,,
\\
v^1_{(1, C)} & = \lran{AC} - x\!\cdot\!\lrbr{C^2}\,,
& \hspace{1ex}
v^1_{(2, C)} & = \lran{BC} - y^{\beta' - \beta}\!\cdot\!\lrbr{AC}\,,
& \hspace{1ex}
v^1_{(3, C)} & = \lran{C^2} - z^{\gamma' - \gamma}\!\cdot\!\lrbr{BC}\,.
\end{alignat*}
We set
$\lsp{\Lambda^{\!1}} =
\{ (1, A)\,,\,(2, A)\,,\,(3, A)\,,\,(2, B)\,,\,(3, B)\,,\,(3, C) \}
\subseteq \widetilde{\Lambda^{\!1}} =
\{ 1, 2, 3 \} \times \{ A, B, C \}$\,.
Then we have the following.

\begin{lem}\label{4.5}
$\{ v^1_{(i, M)} \}_{(i, M)\,\in\,\lsp{\Lambda^{\!1}}}
\,\cup\,
\lrbr{\mon{A, B, C}{2}}$
is an $R$-free basis of $F^1_2$\,.
\end{lem}

\noindent
{\it Proof}.\,
Let us recall that
$\lrbr{\mon{A, B, C}{2}} \cup \lran{\mon{A, B, C}{2}}$
is an $R$-free basis of $F^1_2$\,.
Because $|\lsp{\Lambda^{\!1}}| = |\mon{A, B, C}{2}| = 6$ and
$\lran{\mon{A, B, C}{2}} \subseteq
R\!\cdot\!\{ v^1_{(i, M)} \}_{(i, M)\,\in\,\lsp{\Lambda^{\!1}}}
+ R\!\cdot\!\lrbr{\mon{A, B, C}{2}}$\,,
we get the required assertion by \ref{3.3}.

\vspace{1em}
Let $K_\bullet = K_\bullet(x, y^\beta, z^\gamma)$\,.
By \ref{3.2} there exists a chain map
$\sigma^1_\bullet : K_\bullet \otimes_R F^1_3 \ra F^1_\bullet$
such that
${\rm Im}\,\sigma^1_0 + {\rm Im}\,\varphi^1_1 = I^3 :_R Q$
and
$\sigma^1_2(\check{e}_i \otimes M) = (-1)^i\!\cdot\! v^1_{(i, M)}$
for any $(i, M) \in \widetilde{\Lambda^{\!1}} =
\{ 1, 2, 3 \} \times \{ A, B, C \}$\,.
Moreover, we get an acyclic complex
\[
0 \ra \lsp{F^1_3} \homom{\lsp{\varphi^1_3}}
\lsp{F^1_2} \homom{\lsp{\varphi^1_2}} 
\lss{F^1_1} \homom{\lss{\varphi^1_1}} \lss{F^1_0} = R\,, 
\]
where
\[
\lsp{F^1_3} = K_2 \otimes_R F^1_3\,, \hspace{2ex}
\lsp{F^1_2} = 
\begin{array}{c}
K_1 \otimes_R F^1_3 \\
\oplus \\
F^1_2
\end{array}\,, \hspace{2ex}
\lss{F^1_1} = 
\begin{array}{c}
K_0 \otimes_R F^1_3 \\
\oplus \\
F^1_1
\end{array}\,, \hspace{2ex}
\lsp{\varphi^1_3} = \left(
\begin{array}{c}
\partial_2 \otimes {\rm id}_{F^1_3} \\
\sigma^1_2
\end{array}\right)\,,
\]
and
$\lss{\varphi^1_1} = 
(\,\sigma^1_0 \hspace{1.5ex} \varphi^1_1\,)$\,.
Let us recall our notation introduced in Section 3.
For any $(i, M) \in \widetilde{\Lambda^{\!1}}$ we set
\[
\lrbr{i, M} = \lrbr{e_i \otimes M} =
\binom{e_i \otimes M}{0} \in \lsp{F^1_3}\,.
\]
On the other hand,
for any $\eta \in F^1_2$\,,
we set
\[
\lran{\eta} = \binom{0}{\eta} \in \lsp{F^1_2}\,.
\]
In particular, as $F^1_2 = S_2 \oplus S_2$\,,
$\lran{\lrbr{h}} \in \lsp{F^1_2}$ is defined for any
$h \in S_2$\,.
We set
$\lran{\lrbr{\mon{A, B, C}{2}}} =
\{ \lran{\lrbr{M}} \}_{M \in \mon{A, B, C}{2} }$\,.
Let $\lss{F^1_2}$ be the $R$-submodule of $\lsp{F^1_2}$
generated by
\[
\{ \lrbr{i, M} \}_{(i, M)\,\in\,\widetilde{\Lambda^{\!1}}}
\,\cup\, 
\lran{\lrbr{\mon{A, B, C}{2}}}\,,
\]
and let $\lss{\varphi^1_2}$ be the restriction of
$\lsp{\varphi^1_2}$ to $\lss{F^1_2}$\,.
In order to define $\lss{F^1_3}$\,,
we set $\lss{\Lambda^{\!1}} =
\widetilde{\Lambda^{\!1}} \setminus \lsp{\Lambda^{\!1}} =
\{ (1, B)\,,\,(1, C)\,,\,(2, C) \}$\,.
We need the next result which can be checked directly.

\begin{lem}\label{4.6}
The following equalities hold;
\begin{eqnarray*}
v^1_{(1, B)} & = &
v^1_{(2, A)} - x\!\cdot\!\lrbr{BC} +
y^{\beta' - \beta}\!\cdot\!\lrbr{A^2}\,, \\
v^1_{(1, C)} & = &
v^1_{(3, A)} - x\!\cdot\!\lrbr{C^2} +
z^{\gamma' - \gamma}\!\cdot\!\lrbr{AB}\,, \\
v^1_{(2, C)} & = &
v^1_{(3, B)} - y^{\beta' - \beta}\!\cdot\!\lrbr{AC} +
z^{\gamma' - \gamma}\!\cdot\!\lrbr{B^2}\,.
\end{eqnarray*}
\end{lem}

\noindent
So, we define the elements
$\lss{w^1_{(i, M)}} \in \lsp{F^1_3}$
for $(i, M) \in \lss{\Lambda^{\!1}}$ as follows;
\begin{eqnarray*}
\lss{w^1_{(1, B)}} & = &
-\check{e}_1 \otimes B - \check{e}_2 \otimes A\,, \\
\lss{w^1_{(1, C)}} & = &
-\check{e}_1 \otimes C + \check{e}_3 \otimes A\,, \\
\lss{w^1_{(2, C)}} & = &
\check{e}_2 \otimes C + \check{e}_3 \otimes B\,.
\end{eqnarray*}
Let $\lss{F^1_3}$ be the $R$-submodule of $\lsp{F^1_3}$
generated by
$\{ \lss{w^1_{(i, M)}} \}_{(i, M) \in \lss{\Lambda^{\!1}}}$
and let $\lss{\varphi^1_3}$ be the restriction of
$\lsp{\varphi^1_3}$ to $\lss{F^1_3}$\,.
Thus we get a complex
\[
0 \ra \lss{F^1_3} \homom{\lss{\varphi^1_3}}
\lss{F^1_2} \homom{\lss{\varphi^1_2}} 
\lss{F^1_1} \homom{\lss{\varphi^1_1}} \lss{F^1_0} = R\,.
\]
Let us denote $(\lss{F^1_\bullet}\,, \lss{\varphi^1_\bullet})$ by
$(F^2_\bullet\,, \varphi^2_\bullet)$\,.
Moreover, we put $w^2_{(i, M)} = \lss{w^1_{(i, M)}}$ for
$(i, M) \in \lss{\Lambda^{\!1}}$\,.
Then, by \ref{3.5} and \ref{3.6} we have the next result.

\begin{lem}\label{4.7}
$(F^2_\bullet\,, \varphi^2_\bullet)$
is an acyclic complex of finitely generated free $R$-modules
satisfying the following conditions.
\begin{itemize}
\item[{\rm (1)}]
${\rm Im}\,\varphi^2_1 = I^3 :_R Q$\,.
\item[{\rm (2)}]
$\{ w^2_{(1, B)}\,, w^2_{(1, C)}\,, w^2_{(2, C)} \}$
is an $R$-free basis of $F^2_3$\,.
\item[{\rm (3)}]
$\{ \lrbr{i, M} \}_{(i, M)\,\in\,\widetilde{\Lambda^{\!1}}}
\,\cup\,
\lran{\lrbr{\mon{A, B, C}{2}}}$ is an $R$-free basis of $F^2_2$\,.
\item[{\rm (4)}]
The following equalities hold;
\begin{equation*}
\begin{split}
\varphi^2_3(w^2_{(1, B)}) &= 
- y^\beta\!\cdot\!\lrbr{3, B} + z^\gamma\!\cdot\!\lrbr{2, B} -
x\!\cdot\!\lrbr{3, A} + z^\gamma\!\cdot\!\lrbr{1,A}     \\
  &   \hspace{30ex}
- x\!\cdot\!\lran{\lrbr{BC}} + y^{\beta' - \beta}\!\cdot\!\lran{\lrbr{A^2}}\,, \\
\varphi^2_3(w^2_{(1, C)}) &=
-y^\beta\!\cdot\!\lrbr{3, C} + z^\gamma\!\cdot\!\lrbr{2, C} +
x\!\cdot\!\lrbr{2, A} - y^\beta\!\cdot\!\lrbr{1, A}   \\
  &   \hspace{30ex}
- x\!\cdot\!\lran{\lrbr{C^2}} + z^{\gamma' - \gamma}\!\cdot\!\lran{\lrbr{AB}}\,, \\
\varphi^2_3(w^2_{(2, C)}) &=
x\!\cdot\!\lrbr{3, C} - z^\gamma\!\cdot\!\lrbr{1, C} +
x\!\cdot\!\lrbr{2, B} - y^\beta\!\cdot\!\lrbr{1, B}   \\
  &   \hspace{30ex}
- y^{\beta' - \beta}\!\cdot\!\lran{\lrbr{AC}} + z^{\gamma' - \gamma}\!\cdot\!\lran{\lrbr{B^2}}\,.
\end{split} 
\end{equation*}
\end{itemize}
\end{lem}

We put $\Lambda^{\!2} = \lss{\Lambda^{\!1}} =
\{ (1, B)\,, (1, C)\,, (2, C) \}$ and
$\widetilde{\Lambda^{\!2}} = 
\{ 1, 2, 3 \} \times \Lambda^{\!2}$\,.
We simply denote
$(j, (i, M)) \in \widetilde{\Lambda^{\!2}}$ by
$(j, i, M)$\,. Then
\[
\widetilde{\Lambda^{\!2}} = 
\left\{\begin{array}{lll}
(1, 1, B)\,, & (2, 1, B)\,, & (3, 1, B)\,, \\
(1, 1, C)\,, & (2, 1, C)\,, & (3, 1, C)\,, \\
(1, 2, C)\,, & (2, 2, C)\,, & (3, 2, C)
\end{array}\right\}\,.
\]
As we are assuming
$2\beta \leq \beta'$ and $2\gamma \leq \gamma'$\,,
by (4) of \ref{4.7} we get
\[
\varphi^2_3(w^2_{(i, M)}) =
x\!\cdot\!\varphi^2_{(1,i, M)} +
y^\beta\!\cdot\!\varphi^2_{(2, i, M)} +
z^\gamma\!\cdot\!\varphi^2_{(3, i, M)}
\]
for any $(i, M) \in \Lambda^{\!2}$\,, where
\begin{align*}
 & v^2_{(1, 1, B)}
= -\lrbr{3, A} - \lran{\lrbr{BC}}\,,
\hspace{1.5ex} v^2_{(2, 1, B)}
= -\lrbr{3, B} + y^{\beta'-2\beta}\!\cdot\!\lran{\lrbr{A^2}}\,,
\hspace{1.5ex} v^2_{(3, 1, B)}
= \lrbr{2, B} + \lrbr{1, A}\,, \\
 & v^2_{(1, 1, C)}
= \lrbr{2, A} - \lran{\lrbr{C^2}}\,, 
\hspace{1.5ex} v^2_{(2, 1, C)} 
= -\lrbr{3, C} - \lrbr{1, A}\,, 
\hspace{1.5ex} v^2_{(3, 1, C)} 
= \lrbr{2, C} + z^{\gamma'-2\gamma}\!\cdot\!\lran{\lrbr{AB}}\,, \\
 & v^2_{(1, 2, C)} 
= \lrbr{3, C} + \lrbr{2, B}\,,
\hspace{1.5ex} v^2_{(2, 2, C)} 
= -\lrbr{1, B} - y^{\beta'-2\beta}\!\cdot\!\lran{\lrbr{AC}}
\hspace{1.5ex}\mbox{and}\hspace{1.5ex}\\
 & \hspace{50ex} v^2_{(3, 2, C)} 
= -\lrbr{1, C} + z^{\gamma'-2\gamma}\!\cdot\!\lran{\lrbr{B^2}}\,.
\end{align*}
Thus a family
$\{ v^2_{(j, i, M)} \}_{(j, i, M)\,\in\,\widetilde{\Lambda^{\!2}}}$
of elements in $F^2_2$ is fixed and we see
${\rm Im}\,\varphi^2_3 \subseteq QF^2_2$\,. Because
${\rm Im}\,\varphi^2_1 :_R Q = (I^3 :_R Q) :_R Q = I^3 :_R Q^2$
and $F^2_3 \cong R^{\oplus 3}$\,,
by \ref{3.1} we get the next result.

\begin{thm}\label{4.8}
Let $\alpha = 1$\,, $\alpha' = 2$\,, $2\beta \leq \beta'$ and
$2\gamma \leq \gamma'$\,.
Then we have
\[
(I^3 :_R Q^2) / (I^3 :_R Q) \cong (R / Q)^{\oplus 3}\,.
\]
\end{thm}

The following relation, which can be checked directly,
is very important.

\begin{lem}\label{4.9}
$v^2_{(3, 1, B)} + v^2_{(2, 1, C)} + v^2_{(1, 2, C)}
= 2\!\cdot\!\lrbr{2, B}$\,.
\end{lem}

By \ref{3.2} there exists a chain map
$\sigma^2_\bullet : K_\bullet \otimes_R F^2_3 
\ra F^2_\bullet$ such that
${\rm Im}\,\sigma^2_0 + {\rm Im}\,\varphi^2_1 = I^3 :_R Q^2$
and
$\sigma^2_2(\check{e}_j \otimes w^2_{(i, M)}) =
(-1)^j\!\cdot\! v^2_{(j, i, M)}$
for any $(j, i, M) \in \widetilde{\Lambda^{\!2}}$\,.
Moreover, we get an acyclic complex
\[
0 \ra \lsp{F^2_3} \homom{\lsp{\varphi^2_3}}
\lsp{F^2_2} \homom{\lsp{\varphi^2_2}} 
\lss{F^2_1} \homom{\lss{\varphi^2_1}} \lss{F^2_0} = R\,, 
\]
where
\[
\lsp{F^2_3} = K_2 \otimes_R F^2_3\,, \hspace{2ex}
\lsp{F^2_2} = 
\begin{array}{c}
K_1 \otimes_R F^2_3 \\
\oplus \\
F^2_2
\end{array}\,, \hspace{2ex}
\lss{F^2_1} = 
\begin{array}{c}
K_0 \otimes_R F^2_3 \\
\oplus \\
F^2_1
\end{array}\,, \hspace{2ex}
\lsp{\varphi^2_3} = \left(
\begin{array}{c}
\partial_2 \otimes {\rm id}_{F^2_3} \\
\sigma^2_2
\end{array}\right)\,,
\]
and $\lss{\varphi^2_1} =
(\sigma^2_0 \hspace{1.5ex} \varphi^2_1)$\,.
In order to remove non-minimal components from
$\lsp{F^2_3}$ and $\lsp{F^2_2}$\,,
we would like to choose a subset
$\lsp{\Lambda^{\!2}}$ of $\widetilde{\Lambda^{\!2}}$
as big as possible so that
$\{ v^2_{(j, i, M)} \}_{(j, i, M) \in \lsp{\Lambda^{\!2}}}$
is a part of an $R$-free basis of $F^2_2$\,.

\begin{thm}\label{4.10}
Let $\alpha = 1$\,, $\alpha' = 2$\,, $2\beta \leq \beta'$ and
$2\gamma \leq \gamma'$\,.
Suppose that $2$ is a unit in $R$\,.
Then, we can take $\widetilde{\Lambda^{\!2}}$ itself
as $\lsp{\Lambda^{\!2}}$\,.
Hence
$\dep{}{R / (I^3 :_R Q^2)} > 0$\,,
and so $I^{(3)} = I^3 :_R Q^2$\,.
Moreover, we have
$\length{R}{I^{(3)} / I^3} = 6\!\cdot\!\length{R}{R / Q}$\,.
\end{thm}

\noindent
{\it Proof}.\,
We would like to show that
\[
\{ v^2_{(j, i, M)} \}_{(j, i, M)\,\in\,\widetilde{\Lambda^{\!2}}}
\, \cup \,
 \lran{\lrbr{\mon{A, B, C}{2}}}
\]
is an $R$-free basis of $F^2_2$\,.
Let us recall that
\[
\{ \lrbr{i, M} \}_{(i, M)\,\in\,\widetilde{\Lambda^{\!1}}}
\, \cup \,
 \lran{\lrbr{\mon{A, B, C}{2}}}
\]
is an $R$-free basis of $F^2_2$\,.
Because $|\widetilde{\Lambda^{\!1}} | =
|\widetilde{\Lambda^{\!2}}|$\,,
we have to prove
\[
\{ \lrbr{i, M} \}_{(i, M)\,\in\,\widetilde{\Lambda^{\!1}}}
\, \subseteq  \, G :=
R\!\cdot\!\{ v^2_{(j, i, M)} \}_{(j, i, M)\,\in\,\widetilde{\Lambda^{\!2}}}
+ R\!\cdot\!\lran{\lrbr{\mon{A, B, C}{2}}}\,.
\]
In fact, we have
$\lrbr{2, A} = v^2_{(1, 1, C)} + \lran{\lrbr{C^2}} \in G$\,.
Similarly, we can easily see that
$\lrbr{3, A}$\,, $\lrbr{1, B}$\,, $\lrbr{3, B}$\,, $\lrbr{1, C}$
and $\lrbr{2, C}$ are included in $G$\,.
Moreover, as $2$ is a unit in $R$\,,
we have $\lrbr{2, B} \in G$ by \ref{4.9}.
Then $\lrbr{1, A} = v^2_{(3, 1, B)} - \lrbr{2, B} \in G$
and $\lrbr{3, C} = v^2_{(1, 2, C)} - \lrbr{2, B} \in G$\,.
The last assertion holds since
\begin{eqnarray*}
\length{R}{I^{(3)} / I^3} & = &
 \length{R}{(I^3 :_R Q^2) / I^3} \\
 & = &
 \length{R}{(I^3 :_R Q^2) / (I^3 :_R Q)} +
 \length{R}{(I^3 :_R Q) / I^3} \\
 & = &
 3\!\cdot\!\length{R}{R / Q} + 3\!\cdot\!\length{R}{R / Q}
\end{eqnarray*}
by \ref{4.3} and \ref{4.8}.
Thus the proof is completed.

\vspace{1em}
In the rest of this section,
let us consider the case where ${\rm ch}\,R = 2$\,.
In this case, we have
\[
v^2_{(3, 1, B)} + v^2_{(2, 1, C)} + v^2_{(1, 2, C)} = 0\,.
\]
We set $\lsp{\Lambda^{\!2}} =
\widetilde{\Lambda^{\!2}} \setminus \{\,(3, 1, B)\,\}$\,.
Then, it is easy to see that
\[
\{ v^2_{(j, i, M)} \}_{(j, i, M)\,\in\,\lsp{\Lambda^{\!2}}}
\,\cup\,
\{\,\lrbr{2, B}\,\}
\,\cup\,
\lran{\lrbr{\mon{A, B, C}{2}}}
\]
is an $R$-free basis of $F^2_2$\,.
For any $(j, i, M) \in \widetilde{\Lambda^{\!2}}$\,,
let us simply denote
$\lrbr{j, (i, M)} =
\lrbr{e_j \otimes w^2_{(i, M)}} \in \lsp{F^2_2}$ by
$\lrbr{j, i, M}$\,. Then
\[
\{\, \lrbr{j, i, M} \,\}_{
(j, i, M)\,\in\,\widetilde{\Lambda^{\!2}}}
\,\cup\,
\{\, \lran{v^2_{(j, i, M)}} \,\}_{
(j, i, M)\,\in\,\lsp{\Lambda^{\!2}}}
\,\cup\,
\{\, \lran{\lrbr{2, B}} \,\}
\,\cup\,
\lran{\lran{\lrbr{\mon{A, B, C}{2}}}}
\]
is an $R$-free basis of $\lsp{F^2_2}$\,.
Let $\lss{F^2_2}$ be the $R$-submodule of $\lsp{F^2_2}$
generated by
\[
\{\, \lrbr{j, i, M} \,\}_{
(j, i, M)\,\in\,\widetilde{\Lambda^{\!2}}}
\,\cup\,
\{\, \lran{\lrbr{2, B}} \,\}
\,\cup\,
\lran{\lran{\lrbr{\mon{A, B, C}{2}}}}
\]
and let $\lss{\varphi^2_2}$ be the restriction of
$\lsp{\varphi^2_2}$ to $\lss{F^2_2}$\,.
In order to define $\lss{F^2_3}$\,,
we set
$\lss{\Lambda^{\!2}} =
\widetilde{\Lambda^{\!2}} \setminus \lsp{\Lambda^{\!2}}
= \{ (3, 1, B) \}$\,.
Because
\[
v^2_{(3, 1, B)} = -v^2_{(2, 1, C)} - v^2_{(1, 2, C)}\,,
\]
we define $\lss{w^2_{(3, 1, B)}} \in \lsp{F^2_3}$ to be
\[
-\check{e}_3 \otimes w^2_{(1, B)} +
\check{e}_2 \otimes w^2_{(1, C)} -
\check{e}_1 \otimes w^2_{(2, C)}\,.
\]
Let $\lss{F^2_3}$ be the $R$-submodule of
$\lsp{F^2_3}$ generated by $\lss{w^2_{(3, 1, B)}}$
and let $\lss{\varphi^2_3}$ be the restriction of
$\lsp{\varphi^2_3}$ to $\lss{F^2_3}$\,.
Thus we get a complex
\[
0 \ra \lss{F^2_3} \homom{\lss{\varphi^2_3}}
\lss{F^2_2} \homom{\lss{\varphi^2_2}} 
\lss{F^2_1} \homom{\lss{\varphi^2_1}} \lss{F^2_0} = R\,.
\]
Let us denote
$(\lss{F^2_\bullet}\,, \lss{\varphi^2_\bullet})$ by
$(F^3_\bullet\,, \varphi^3_\bullet)$\,.
Moreover, we put
$w^3_{(3, 1, B)} = \lss{w^2_{(3, 1, B)}}$\,.
Then, by \ref{3.5} and \ref{3.6} we have the next result.

\begin{lem}\label{4.11}
$(F^3_\bullet\,, \varphi^3_\bullet)$
is an acyclic complex of finitely generated free $R$-modules
satisfying the following conditions.
\begin{itemize}
\item[{\rm (1)}]
${\rm Im}\,\varphi^3_1 = I^3 :_R Q^2$\,.
\item[{\rm (2)}]
$w^3_{(3, 1, B)}$ is an $R$-free basis of $F^3_3$\,.
\item[{\rm (3)}]
$\{\, \lrbr{j, i, M} \,\}_{
(j, i, M)\,\in\,\widetilde{\Lambda^{\!2}}}
\,\cup\,
\{\,\lran{\lrbr{2, B}}\,\}
\,\cup\,
\lran{\lran{\lrbr{\mon{A, B, C}{2}}}}$
is an $R$-free basis of $F^3_2$\,.
\item[{\rm (4)}]
The following equality holds {\rm ;}
\begin{equation*}
\begin{split}
\varphi^3_3(w^3_{(3, 1, B)}) &=
- x\!\cdot\!\lrbr{2, 1, B} + y^\beta\!\cdot\!\lrbr{1, 1, B} 
+ x\!\cdot\!\lrbr{3, 1, C}  \\
 & \hspace{10ex} 
- z^\gamma\!\cdot\!\lrbr{1, 1, C} 
- y^\beta\!\cdot\!\lrbr{3, 2, C} + z^\gamma\!\cdot\!\lrbr{2, 2, C}\,.
\end{split}
\end{equation*}
\end{itemize}
\end{lem} 

We put $\Lambda^{\!3} = \lss{\Lambda^{\!2}} =
\{ (3, 1, B) \}$ and
$\widetilde{\Lambda^{\!3}} =
\{ 1, 2, 3 \} \times \Lambda^{\!3}$\,.
We simply denote $(i, (3, 1, B)) \in
\widetilde{\Lambda^{\!3}}$
by $(i, 3, 1, B)$\,.
Then $\widetilde{\Lambda^{\!3}} =
\{\,(i, 3, 1, B)\,\}_{i = 1, 2, 3}$\,.
By (4) of \ref{4.11} we have
\[
\varphi^3_3(w^3_{(3, 1, B)}) =
x\!\cdot\! v^3_{(1, 3, 1, B)} +
y^\beta\!\cdot\! v^3_{(2, 3, 1, B)} +
z^\gamma\!\cdot\! v^3_{(3, 3, 1, B)}\,,
\]
where
\begin{equation*}
\begin{split}
 & v^3_{(1, 3, 1, B)} =
\lrbr{3, 1, C} - \lrbr{2, 1, B} \,,
\hspace{2ex}
v^3_{(2, 3, 1, B)} = 
\lrbr{1, 1, B} - \lrbr{3, 2, C} 
\hspace{2ex}\mbox{and}\hspace{2ex} \\
 & \hspace{45ex}
v^3_{(3, 3, 1, B)} =
\lrbr{2, 2, C} - \lrbr{1, 1, C}\,.
\end{split}
\end{equation*}
Thus a family
$\{ v^3_{(i, 3, 1, B)} \}_{i = 1, 2, 3}$
of elements in $F^3_2$
is fixed and we see ${\rm Im}\,\varphi^3_3 \subseteq QF^3_2$\,.
Because
\[
{\rm Im}\,\varphi^3_1 :_R Q = (I^3 :_R Q^2) :_R Q =
I^3 :_R Q^3
\]
and $F^3_3 \cong R$\,,
by \ref{3.1} we get the next result.

\begin{thm}\label{4.12}
Let $\alpha = 1$\,, $\alpha' = 2$\,, $2\beta \leq \beta'$\,,
$2\gamma \leq \gamma'$ and ${\rm ch}\,R = 2$\,.
Then we have
\[
(I^3 :_R Q^3) / (I^3 :_R Q^2) \cong R / Q\,.
\]
\end{thm}

It is easy to see that
\[
\{ v^2_{(i, 3, 1, B)} \}_{i = 1, 2, 3} 
\,\cup\,
\left\{\begin{array}{ll}
\lrbr{2, 1, B}\,, & \lrbr{3, 1, B}\,, \\
\lrbr{1, 1, C}\,, & \lrbr{2, 1, C}\,, \\
\lrbr{1, 2, C}\,, & \lrbr{3, 2, C}
\end{array}\right\}
\,\cup\,
\lran{\lran{\lrbr{\mon{A, B, C}{2}}}}
\]
is an $R$-free basis of $F^3_2$\,.
Therefore, by \ref{3.4} we get the following result.

\begin{thm}\label{4.13}
Let $\alpha = 1$\,, $\alpha' = 2$\,, $2\beta \leq \beta'$\,,
$2\gamma \leq \gamma'$ and ${\rm ch}\,R = 2$\,.
Then, it follows that
$\dep{}{R / (I^3 :_R Q^3)} > 0$\,,
and so $I^{(3)} = I^3 :_R Q^3$\,.
Moreover, we have
$\length{R}{I^{(3)} / I^3} = 7\!\cdot\!\length{R}{R / Q}
= 7\alpha\beta\!\cdot\!\length{R}{R / (x, y, z)R}$\,.
\end{thm}

\noindent
The last assertion of \ref{4.13} holds since $\length{R}{I^{(3)} / I^3}$ coincides with
\begin{eqnarray*}
 &   & \length{R}{(I^3 :_R Q^3) / (I^3 :_R Q^2)} +
 \length{R}{(I^3 :_R Q^2) / (I^3 :_R Q)} +
 \length{R}{(I^3 :_R Q) / I^3}   \\
 &  = &
 \length{R}{R / Q} + 3\!\cdot\!\length{R}{R / Q} + 3\!\cdot\!\length{R}{R / Q}
\end{eqnarray*}
by \ref{4.3}, \ref{4.8} and \ref{4.12}.

\section{Computing $\epsilon$-multiplicity}

Let $R$ be a $3$-dimensional regular local ring
with the maximal ideal $\gm = (x, y, z)R$\,.
Let $I$ be the ideal generated by the maximal minors
of the matrix
\[
\left(\begin{array}{lll}
x & y & z \\
y & z & x^2
\end{array}\right)\,.
\]
We will compute the length of $I^{(n)} / I^n$ for all $n$
using the $\ast$-transforms.
As a consequence of our result, we get $\epsilon(I) = 1 / 2$\,, where
\[
\epsilon(I) := \lim_{n \rightarrow \infty} \frac{3!}{n^3}\!\cdot\!
\length{R}{I^{(n)} / I^n}\,,
\]
which is the invariant called $\epsilon$-multiplicity of $I$
(cf. \cite{chs}, \cite{uv}).
Let us maintain the same notations as in Section 4.
So, $a = z^2 - x^2y$\,, $b = x^3 - yz$\,, $c = y^2 - xz$\,,
$S = R[A, B, C]$\,, $f = xA + yB + zC$ and $g = yA + zB + x^2C$\,.
Furthermore, we need the following notations which is not used in Section 4.
For any $0 \leq d \in \zz$\,, we denote the set $\{ B^\beta C^\gamma 
\mid 0 \leq \beta\,, \gamma \in \zz
\hspace{1ex}\mbox{and}\hspace{1ex}
\beta + \gamma = d \}$
by $\mon{B, C}{d}$\,,
and for a monomial $L$ and a set ${\cal S}$ consisting of monomials,
we denote the set
$\{ LM \mid M \in {\cal S} \}$ by $L\!\cdot\!{\cal S}$\,.

Let $n \geq 2$\,.
Then the complex
\[
0 \ra S_{n-2} \homom{\binom{-g}{f}} S_{n-1} \oplus S_{n-1}
\homom{(f \hspace{0.5ex} g)} S_n \homom{\epsilon} R
\]
is acyclic and gives a minimal free resolution of $I^n$\,.
We denote this complex by $(F^1_\bullet\,, \varphi^1_\bullet)$\,.

As $\Lambda^{\!1}$, which is an $R$-free basis of $F^1_3 = S_{n-2}$\,,
we take $\mon{A, B, C}{n-2}$\,.
Then, for any $M \in \Lambda^{\!1}$\,, we have
\[
\varphi^1_3(M) = x\!\cdot\! v^1_{(1, M)} +
y\!\cdot\! v^1_{(2, M)} + z\!\cdot\! v^1_{(3, M)}\,,
\]
where
\[
v^1_{(1, M)} = \lran{AM} - x\!\cdot\!\lrbr{CM}\,,\hspace{1ex}
v^1_{(2, M)} = \lran{BM} - \lrbr{AM}\,,\hspace{1ex}
v^1_{(3, M)} = \lran{CM} - \lrbr{BM}\,.
\]
Hence, by \ref{3.1} we get the following.

\begin{lem}\label{5.1}
$(I^n :_R \gm) / I^n \cong (R / \gm)^{\oplus\binom{n}{2}}$ \,.
\end{lem}

The next result can be checked directly.

\begin{lem}\label{5.2}
Suppose $n \geq 4$\,.
Then, for any $N \in \mon{A, B, C}{n-4}$\,, we have
\[
v^1_{(3, A^2N)} = v^1_{(2, ABN)} - v^1_{(1, B^2N)} + v^1_{(1, ACN)} -
x\!\cdot\! v^1_{(2, C^2N)} + x\!\cdot\! v^1_{(3, BCN)}\,.
\]
\end{lem}

Let $\widetilde{\Lambda^{\!1}} =
\{ 1, 2, 3 \} \times \Lambda^{\!1}$\,.
Then the following result holds.

\begin{lem}\label{5.3}
As an $R$-module, 
$F^1_2 = S_{n-1} \oplus S_{n-1}$ is generated by
\[
\{ v^1_{(i, M)} \}_{
(i, M)\,\in\,\widetilde{\Lambda^{\!1}}}
\,\cup\,
\{ \lrbr{AC^{n-2}}\,, \lrbr{BC^{n-2}}\,,
\lrbr{C^{n-1}} \}\,.
\]
\end{lem}

\noindent
{\it Proof}.\,
Let $G$ be the sum of the
$R$-submodule of $F^1_2$ generated by
the elements stated above and $\gm F^1_2$.
It is enough to show that
$\lrbr{\mon{A, B, C}{n-1}}$ and
$\lran{\mon{A, B, C}{n-1}}$ are contained in
$G$\,.

First, let us prove
$\lran{L} \in G$
for any $L \in \mon{A, B, C}{n-1}$\,.
We write $L = A^\alpha B^\beta C^\gamma$\,,
where $0 \leq \alpha, \beta, \gamma \in \zz$ and
$\alpha + \beta + \gamma = n - 1$\,.
If $\alpha > 0$\,,
\[
\lran{L} = \lran{A\!\cdot\! A^{\alpha - 1}B^\beta C^\gamma}
= v^1_{(1, A^{\alpha - 1}B^\beta C^\gamma)}
+ x\!\cdot\!\lrbr{C\!\cdot\! A^{\alpha - 1}B^\beta C^\gamma}
\in G\,.
\]
Hence, we have to consider the case where $\alpha = 0$\,.
However, as
\begin{equation*}
\begin{split}
 & \lran{C^{n-1}} = \lran{C\!\cdot\! C^{n-2}} =
v^1_{(3, C^{n-2})} + \lrbr{B\!\cdot\! C^{n-2}} \in G
\hspace{2ex}\mbox{and} \\
 & \lran{BC^{n-2}} = \lran{B\!\cdot\! C^{n-2}} =
v^1_{(2, C^{n-2})} + \lrbr{A\!\cdot\! C^{n-2}} \in G\,,
\end{split}
\end{equation*}
we may assume $\beta \geq 2$\,.
Then, as
\begin{eqnarray*}
\lran{L} & = & \lran{B\!\cdot\! B^{\beta - 1}C^\gamma} \\
 & = & v^1_{(2, B^{\beta - 1}C^\gamma)} +
       \lrbr{A\!\cdot\! B^{\beta - 1}C^\gamma}  \\
 & = & v^1_{(2, B^{\beta - 1}C^\gamma)} +
       \lrbr{B\!\cdot\! AB^{\beta - 2}C^\gamma}  \\
 & = & v^1_{(2, B^{\beta - 1}C^\gamma)} -
       v^1_{(3, AB^{\beta - 2}C^\gamma)} +
       \lran{C\!\cdot\! AB^{\beta - 2}C^\gamma}
\end{eqnarray*}
and as
$\lran{C\!\cdot\! AB^{\beta - 2}C^\gamma} =
\lran{A\!\cdot\! B^{\beta - 2}C^{\gamma + 1}}$\,,
we get $\lran{L} \in G$\,.

Next, we prove
$\lrbr{\mon{A, B, C}{n-1}} \subseteq G$\,.
Let us notice
$\mon{A, B, C}{n-1} =
A\!\cdot\!\mon{A, B, C}{n-2} \cup B\!\cdot\!\mon{B, C}{n-2}
\cup \{ C^{n-1} \}$\,.
For any $M \in \mon{A, B, C}{n-2}$ and any
$X \in \mon{B, C}{n-2}$\,, we have
\[
\lrbr{AM} = -v^1_{(2, M)} + \lran{BM} \in G
\hspace{3ex}\mbox{and}\hspace{3ex}
\lrbr{BX} = -v^1_{(3, X)} + \lran{CX} \in G\,.
\]
Hence, the proof is complete as
$\lrbr{C^{n-1}} \in G$ holds obviously.

\vspace{1em}
Now, let $q$ be the largest integer such that $q \leq n / 2$\,.
For any $1 \leq k \leq q$\,,
we would like to construct an acyclic complex
\[
0 \ra F^k_3 \homom{\varphi^k_3} F^k_2
\homom{\varphi^k_2} F^k_1 \homom{\varphi^k_1}
F^k_0 = R
\]
of finitely generated free $R$-modules satisfying
the following conditions.
\begin{itemize}
\item[($\sharp^k_1$)]
${\rm Im}\,{\varphi^k_1} = I^n : \gm^{k - 1}$\,.  
\item[($\sharp^k_2$)]
$F^k_3$ has an $R$-free basis indexed by
$\Lambda^{\!k} := {\rm m}_{A, B, C}^{n - 2k}$\,,
say $\{ w^k_M \}_{M\,\in\,\Lambda^{\!k}}$\,. 
\item[($\sharp^k_3$)]
Let $\widetilde{\Lambda^{\!k}} = \{ 1, 2, 3 \} \times \Lambda^{\!k}$\,.
Then, there exists a family 
$\{ v^k_{(i, M)} \}_{(i, M)\,\in\,\widetilde{\Lambda^{\!k}}}$
of elements in $F^k_2$ satisfying the following conditions\,.
\begin{itemize}
\item[(i)]
For any $M \in \Lambda^{\!k}$\,,
$\varphi^k_3(w^k_{(i, M)}) =
x\!\cdot\! v^k_{(1, M)} + y\!\cdot\! v^k_{(2, M)} + z\!\cdot\! v^k_{(3, M)}$\,.
\item[(ii)]
If $k < q$\,,
for any $N \in \Lambda^{\!{k+1}} := \mon{A, B, C}{n-2k-2}$\,,
\[
v^k_{(3, A^2N)} = v^k_{(2, ABN)} - v^k_{(1, B^2N)} +
v^k_{(1, ACN)} - x\!\cdot\! v^k_{(2, C^2N)} + x\!\cdot\! v^k_{(3, BCN)}\,.
\]
\item[(iii)]
There exists a subset $U^k$ of $F^k_2$ such that
$\{ v^k_{(i, M)} \}_{(i, M)\,\in\,\widetilde{\Lambda^{\!k}}}
\,\cup\, U^k$ generates $F^k_2$ and
\[
| U^k | = 
{\rm rank}\,F^k_2 - 3\!\cdot\!\binom{n - 2k + 2}{2} + \binom{n - 2k}{2}\,,
\]
where the last binomial coefficient is regarded as $0$ if $k = q$\,.
\end{itemize}
\end{itemize}

Let us notice that the acyclic complex
$(F^1_\bullet\,, \varphi^1_\bullet)$\,,
which is already constructed, satisfies 
($\sharp^1_1$), ($\sharp^1_2$) 
($w^1_M$ is $M$ itself
for $M \in \Lambda^{\!1}$)
and ($\sharp^1_3$). 
So, we assume $1 \leq k < q$ 
and an acyclic complex
$(F^k_\bullet\,, \varphi^k_\bullet)$
satisfying the required conditions is given.
Taking the $\ast$-transform of
$(F^k_\bullet\,, \varphi^k_\bullet)$
with respect to $x, y, z$\,,
we would like to construct
$(F^{k+1}_\bullet\,, \varphi^{k+1}_\bullet)$\,.

First, we have the following result since the conditions
($\sharp^k_1$), ($\sharp^k_2$) and (i) of ($\sharp^k_3$)
are satisfied and
$(I^n :_R \gm^{k-1}) :_R \gm = I^n :_R \gm^k$\,.

\begin{lem}\label{5.4}
$(I^n :_R \gm^k) / (I^n :_R \gm^{k-1}) \cong
F_3^k / \gm F^k_3$\,, so
\[
\length{R}{(I^n :_R \gm^k) / (I^n :_R \gm^{k-1})} =
\binom{n - 2k + 2}{2}\,.
\]
\end{lem}

If $\Gamma$ is a subset of $\Lambda^{\!k}$ and $1 \leq i \leq 3$\,,
we denote by $(i, \Gamma)$ the subset
$\{ (i, M) \mid M \in \Gamma \}$ of $\widetilde{\Lambda^{\!k}}$\,.
Let us notice that $\Lambda^{\!k}$ is a disjoint union of
$A^2\!\cdot\!\Lambda^{\!{k+1}}$, $A\!\cdot\!\mon{B, C}{n-2k-1}$ and
$\mon{B, C}{n-2k}$\,.
We set
\[
\lsp{\Lambda^{\!k}} = (1, \Lambda^{\!k}) \,\cup\,
(2, \Lambda^{\!k}) \,\cup\,
(3, A\!\cdot\!\mon{B, C}{n-2k-1} \cup \mon{B, C}{n-2k})\,.
\]
Then the next result holds.

\begin{lem}\label{5.5}
$\{ v^k_{(i, M)} \}_{(i, M)\,\in\,\lsp{\Lambda^{\!k}}} \,\cup\, U^k$
is an $R$-free basis of $F^k_2$\,.
\end{lem}

\noindent
{\it Proof}.\,
Because
\begin{eqnarray*}
| \lsp{\Lambda^{\!k}} | & = &
2\cdot | \Lambda^{\!k} | +
| A\!\cdot\!\mon{B, C}{n-2k-1}\,\cup\,
\mon{B, C}{n-2k} |   \\
 & = & 2\cdot | \Lambda^{\!k} | +
| \Lambda^{\!k} \setminus A^2\!\cdot\!\Lambda^{\!{k+1}} | \\
 & = & 3\cdot | \Lambda^{\!k} | - | \Lambda^{\!{k+1}} | \\
 & = & 3\cdot\binom{n-2k+2}{2} - \binom{n-2k}{2}\,,
\end{eqnarray*}
by (iii) of ($\sharp^k_3$) we have
$| \lsp{\Lambda^{\!k}} | + | U^k | = {\rm rank}\,F^k_2$\,.
Hence, by \ref{3.3} it is enough to show that,
for any $N \in \Lambda^{\!{k+1}}$\,,
$v^k_{(3, A^2N)}$ is contained in the sum of the
$R$-submodule of $F^k_3$
generated by $\{ v^k_{(i, M)} \}_{
(i, M)\,\in\,\lsp{\Lambda^{\!k}}}\,\cup\,U^k$ and
$\gm F^k_2$\,.
We write $N = A^\alpha X$\,,
where $X \in \mon{B, C}{n-2k-2-\alpha}$\,.
Then, using the equalities in (ii) of ($\sharp^k_3$),
the required containment can be proved by induction on $\alpha$\,.

\vspace{1em}
Let $\lss{F^k_2}$ be the $R$-submodule of $\lsp{F^k_2}$
generated by
\[
\{ [i, M] \}_{
(i, M)\,\in\,\widetilde{\Lambda^{\!k}}} \,\cup\, \lran{U^k}\,,
\]
where $[i, M] = [ e_i \otimes w^k_M ]$ for any
$(i, M) \in \widetilde{\Lambda^{\!k}}$\,,
and let $\lss{\varphi^k_2}$ be the restriction of $\lsp{\varphi^k_2}$
to $\lss{F^k_2}$\,.
In order to define $\lss{F^k_3}$\,,
we notice $\widetilde{\Lambda^{\!k}} \setminus \lsp{\Lambda^{\!k}}
= \{ (3, A^2N) \mid N \in \Lambda^{\!{k+1}} \}$\,.
Looking at (ii) of ($\sharp^k_3$), we define
$\lss{w^k_{(3, A^2N)}} \in \lsp{F^k_3}$ to be
\[
-\check{e}_3 \otimes w^k_{A^2N} - \check{e}_2 \otimes w^k_{ABN} -
\check{e}_1 \otimes w^k_{B^2N} + \check{e}_1 \otimes w^k_{ACN} +
x\cdot\check{e}_2 \otimes w^k_{C^2N} + x\cdot\check{e}_3 \otimes w^k_{BCN}
\]
for any $N \in \Lambda^{\!{k+1}}$\,.
Let $\lss{F^k_3}$ be the $R$-submodule of $\lsp{F^k_3}$
generated by $\{ \lss{w^k_{(3, A^2N)}} \}_{
N\,\in\,\Lambda^{\!{k+1}}}$
and let $\lss{\varphi^k_3}$ be the restriction of $\lsp{\varphi^k_3}$
to $\lss{F^k_3}$\,.
Thus we get a complex
\[
0 \ra \lss{F^k_3} \homom{\lss{\varphi^k_3}} \lss{F^k_2}
\homom{\lss{\varphi^k_2}} \lss{F^k_1} \homom{\lss{\varphi^k_1}}
\lss{F^k_0} = R\,.
\]
Let us denote $(\lss{F^k_\bullet}\,, \lss{\varphi^k_\bullet})$ by
$(F^{k+1}_\bullet\,, \varphi^{k+1}_\bullet)$\,.
Moreover, we put $w^{k+1}_N = \lss{w^k_{(3, A^2N)}}$
for any $N \in \Lambda^{\!{k+1}}$\,.
Then, by \ref{3.5} and \ref{3.6} we have the next result.

\begin{lem}\label{5.6}
$(F^{k+1}_\bullet\,, \varphi^{k+1}_\bullet)$
is an acyclic complex satisfying the following conditions.
\begin{itemize}
\item[{\rm (1)}]
${\rm Im}\,\varphi^{k+1}_1 = I^n :_R \gm^k$\,.
\item[{\rm (2)}]
$\{ w^{k+1}_N \}_{N\,\in\,\Lambda^{\!{k+1}}}$
is an $R$-free basis of $F^{k+1}_3$\,.
\item[{\rm (3)}]
$\{ \lrbr{i, M} \}_{(i, M)\,\in\,\widetilde{\Lambda^{\!k}}}
\,\cup\, \lran{U^k}$
is an $R$-free basis of $F^{k+1}_2$\,.
\item[{\rm (4)}]
For any $N \in \Lambda^{\!{k+1}}$\,,
the following equality holds ;
\begin{equation*}
\begin{split}
\varphi^{K+1}_3(w^{k+1}_N) &=
  - x\cdot\lrbr{2, A^2N} + y\cdot\lrbr{1, A^2N}
  - x\cdot\lrbr{3, ABN} + z\cdot\lrbr{1, ABN}  \\
 &\hspace{5ex}
  - y\cdot\lrbr{3, B^2N} + z\cdot\lrbr{2, B^2N}
  + y\cdot\lrbr{3, ACN} - z\cdot\lrbr{2, ACN}  \\
 &\hspace{5ex}
  + x^2\lrbr{3, C^2N} - xz\cdot\lrbr{1, C^2N}
  + x^2\lrbr{2, BCN} - xy\cdot\lrbr{1, BCN}\,.
\end{split}
\end{equation*}
\end{itemize}
\end{lem}

\noindent
The assertions (1) and (2) of the lemma above imply that
$(F^{k+1}_\bullet\,, \varphi^{k+1}_\bullet)$
satisfies ($\sharp^{k+1}_1$) and ($\sharp^{k+1}_2$)\,, respectively.
Moreover, by (4) we have
\[
\varphi^{k+1}_3(w^{k+1}_N) = 
x\cdot v^{k+1}_{(1, N)} + y\cdot v^{k+1}_{(2, N)}
+ z\cdot v^{k+1}_{(3, N)}
\]
for any $N \in \Lambda^{\!{k+1}}$, where
\begin{eqnarray*}
v^{k+1}_{(1, N)} & = &
-\lrbr{2, A^2N} - \lrbr{3, ABN} + 
x\cdot\lrbr{3, C^2N} + x\cdot\lrbr{2, BCN}\,, \\
v^{k+1}_{(2, N)} & = &
\lrbr{1, A^2N} - \lrbr{3, B^2N} +
\lrbr{3, ACN} - x\cdot\lrbr{1, BCN}\,, \\
v^{k+1}_{(3, N)} & = &
\lrbr{1, ABN} + \lrbr{2, B^2N} -
\lrbr{2, ACN} - x\cdot\lrbr{1, C^2N}\,.
\end{eqnarray*}
Thus a family
$\{ v^{k+1}_{(i, N)} \}_{
(i, N)\,\in\,\widetilde{\Lambda^{\!{k+1}}}}$
of elements in $F^{k+1}_2$ satisfying (i) of
($\sharp^{k+1}_3$) is fixed,
where $\widetilde{\Lambda^{\!{k+1}}} =
\{ 1, 2, 3 \} \times \Lambda^{\!{k+1}}$\,.
The next result, which can be checked directly,
insists that (ii) of ($\sharp^{k+1}_3$)
is satisfied if $k + 1 < q$\,.

\begin{lem}\label{5.7}
Suppose $k + 1 < q$\,.
Then $n - 2k - 4 \geq 0$ and we have
\[
v^{k+1}_{(3, A^2L)} =
v^{k+1}_{(2, ABL)} - v^{k+1}_{(1, B^2L)} + v^{k+1}_{(1, ACL)}
- x\cdot v^{k+1}_{(2, C^2L)} + x\cdot v^{k+1}_{(3, BCL)}
\]
for any $L \in \Lambda^{\!{k+2}} := \mon{A, B, C}{n-2k-4}$\,.
\end{lem}

If $\Gamma$ is a subset of $\Lambda^{\!k}$ and $1 \leq i \leq 3$\,,
we denote by $\lrbr{\,i, \Gamma\,}$
the family $\{ \lrbr{i, M} \}_{M\,\in\,\Gamma}$
of elements in $K_1 \otimes_R F^k_3$\,.
The next result means that (iii) of ($\sharp^{k+1}_3$) is satisfied.

\begin{lem}\label{5.8}
We set
\begin{equation*}
\begin{split}
U^{k+1} &= \lrbr{1, A\!\cdot\!\mon{B, C}{n-2k-1}
\,\cup\, \mon{B, C}{n-2k}} \,\cup\,
\lrbr{\,3, \Lambda^{\!k}\,} \,\,\cup\, \\
 &\hspace{5ex}
 \{ \lrbr{2, ABC^{n-2k-2}}\,, \lrbr{2, AC^{n-2k-1}}\,,
 \lrbr{2, BC^{n-2k-1}}\,, \lrbr{2, C^{n-2k}} \} \,\cup\,
 \lran{U^k}\,.
\end{split}
\end{equation*}
Then $\{ v^{k+1}_{(i, N)} \}_{
(i, N)\,\in\,\widetilde{\Lambda^{\!{k+1}}}} \,\cup\,
U^{k+1}$ generates $F^{k+1}_2$ and
\[
| U^{k+1} | = {\rm rank}\,F^{k+1}_2 - 3\!\cdot\!\binom{n - 2k}{2}
+ \binom{n - 2k - 2}{2}\,.
\]
\end{lem}

\noindent
{\it Proof}.\,
Let $G$ be the sum of the $R$-submodule of $F^{k+1}_2$
generated by
$\{ v^{k+1}_{(i, N)} \}_{
(i, N)\,\in\,\widetilde{\Lambda^{\!{k+1}}}} \,\cup\,
U^{k+1}$ and $\gm F^{k+1}_2$\,.
We would like to show $G = F^{k+1}_2$\,.
Let us recall that
\[
\lrbr{\,1, \Lambda^{\!k}\,} \,\cup\,
\lrbr{\,2, \Lambda^{\!k}\,} \,\cup\,
\lrbr{\,3, \Lambda^{\!k}\,} \,\cup\,
\lran{U^k}
\]
is an $R$-free basis of $F^{k+1}_2$ and notice that
$\Lambda^{\!k}$ is a disjoint union of
$A^2\!\cdot\!\Lambda^{\!{k+1}}$\,,
$A\!\cdot\!\mon{A, B, C}{n-2k-1}$ and
$\mon{B, C}{n-2k}$\,.
Because $\lrbr{3, \Lambda^{\!k}} \subseteq U^{k+1}$\,,
it is enough to show
$\lrbr{\,1, A^2\!\cdot\!\Lambda^{\!{k+1}}\,} \,\cup\,
\lrbr{\,2, \Lambda^{\!k}\,} \subseteq G$\,.

First, we prove
$\lrbr{\,1, A^2\!\cdot\!\Lambda^{\!{k+1}}\,} \subseteq G$\,.
Let us take any $N \in \Lambda^{\!{k+1}}$\,.
Then
\[
\lrbr{1, A^2N} = v^{k+1}_{(2, N)} + \lrbr{3, B^2N} -
\lrbr{3, ACN} + x\cdot\lrbr{1, BCN} \in G\,,
\]
and so the required inclusion follows.

Next, we prove $\lrbr{\,2, \Lambda^{\!k}} \subseteq G$\,.
Because
\[
\lrbr{2, A^2N} = -v^{k+1}_{(1, N)} - \lrbr{3, ABN} +
x\cdot\lrbr{3, C^2N} + x\cdot\lrbr{2, BCN} \in G
\]
for any $N \in \Lambda^{\!{k+1}}$\,,
we have $\lrbr{2, A^2\!\cdot\!\Lambda^{\!{k+1}}} \subseteq G$\,.
Furthermore, for any
$B^\beta C^\gamma \in \mon{B, C}{n-2k-1}$\,,
we get $\lrbr{2, AB^\beta C^\gamma} \in G$\,.
In fact,
$\lrbr{2, AB^\beta C^\gamma} \in U^{k+1}$ if
$\beta = \mbox{$0$ or $1$}$\,,
and if $\beta \geq 2$\,, we have
\begin{equation*}
\begin{split}
\lrbr{2, AB^\beta C^\gamma} &=
\lrbr{2, B^2\!\cdot\!AB^{\beta-2}C^\gamma}  \\
 &=
v^{k+1}_{(3, AB^{\beta-2}C^\gamma)} -
\lrbr{1, A^2B^{\beta-1}C^\gamma} +
\lrbr{2, A^2B^{\beta-2}C^{\gamma+1}} +  \\
 &\hspace{35ex}
x\!\cdot\!\lrbr{1, AB^{\beta-2}C^{\gamma+2}} \in G\,.
\end{split}
\end{equation*}
Hence $\lrbr{2, A\!\cdot\!\mon{B, C}{n-2k-1}} \subseteq G$\,.
In order to prove
$\lrbr{2, \mon{B, C}{n-2k}} \subseteq G$\,,
we newly take any
$B^\beta C^\gamma \in \mon{B, C}{n-2k}$\,.
We have $\lrbr{2, B^\beta C^\gamma} \in U^{k+1}$ if
$\beta = \mbox{$0$ or $1$}$\,,
and if $\beta \geq 2$\,, we have
\begin{eqnarray*}
\lrbr{2, B^\beta C^\gamma} & = & 
\lrbr{2, B^2\!\cdot\!B^{\beta-2}C^\gamma}  \\
 & = &
v^{k+1}_{(3, B^{\beta-2}C^\gamma)} - 
\lrbr{1, AB^{\beta-1}C^\gamma} +
\lrbr{2, AB^{\beta-2}C^{\gamma+1}} +
x\!\cdot\!\lrbr{1, B^{\beta-2}C^{\gamma+2}} \in G\,.
\end{eqnarray*}
Hence the required inclusion follows,
and we have seen the first assertion of the theorem.

By (3) of \ref{5.6} we have
${\rm rank}\,F^{k+1}_2 =
3\!\cdot\! | \Lambda^{\!k} | + | U^k |$\,.
On the other hand,
\begin{eqnarray*}
| U^{k+1} | & = &
| \Lambda^{\!k} \setminus A^2\!\cdot\!\Lambda^{\!{k+1}} | +
| \Lambda^{\!k} | + 4 + | U^k | \\
 & = &
2\!\cdot\! | \Lambda^{\!k} | -
| \Lambda^{\!{k+1}} | + 4 + | U^k |\,.
\end{eqnarray*}
Hence we get
\begin{eqnarray*}
{\rm rank}\,F^{k+1}_2 - | U^{k + 1} | & = &
| \Lambda^{\!k} | + | \Lambda^{\!{k+1}} | - 4 \\
 & = &
\binom{n -2k  + 2}{2} + \binom{n - 2k}{2} - 4 \\
 & = &
3\!\cdot\!\binom{n -2k}{2} - \binom{n -2k - 2}{2}\,,
\end{eqnarray*}
and so the second assertion holds.

\vspace{1em}
Thus we have constructed an acyclic complex
\[
0 \ra F^q_3 \homom{\varphi^q_3} F^q_2
\homom{\varphi^q_2} F^q_1 \homom{\varphi^q_1}
F^q_0 = R
\]
of finitely generated free $R$-modules satisfying
($\sharp^q_1$), ($\sharp^q_2$) and ($\sharp^q_3$).
Of course, $n - 2q = \mbox{$0$ or $1$}$\,,
and
\[
\Lambda^{\!q} = \left\{
\begin{array}{ll}
\{\,1\,\} & \mbox{if $n - 2q = 0$} \\
\{ A, B, C \} & \mbox{if $n - 2q = 1$}\,.
\end{array}
\right.
\]
The second condition of (iii) of ($\sharp^q_3$) implies
${\rm rank}\,F^q_2 = 
| \widetilde{\Lambda^{\!q}} | + | U^q |$\,.
Hence, by the first condition of (iii) of ($\sharp^q_3$),
we see that
$\{ v^q_{(i, M)} \}_{
(i, M)\,\in\,\widetilde{\Lambda^{\!q}}} \,\cup\,
\lran{U^q}$
must be an $R$-free basis of $F^q_2$\,.
Therefore, by \ref{3.4} we get the next result.

\begin{thm}\label{5.9}
$\dep{}{R / (I^n :_R \gm^q)} > 0$\,,
and so $I^{(n)} = I^n :_R \gm^q$\,.
\end{thm}

Let us compute $\length{R}{I^{(n)} / I^n}$\,. 
By \ref{5.9} and \ref{5.4} we have
\[
\length{R}{I^{(n)} / I^n} =
\sum_{k = 1}^q \, \length{R}{(I^n : \gm^k) / (I^n : \gm^{k - 1})} =
\sum_{k = 1}^q \, \binom{ n - 2k + 2 }{2}\,.
\]
As a consequence, we get the next result.

\begin{thm}\label{5.10}
The following equality holds ;
\[
\length{R}{I^{(n)} / I^n} = \left\{
\begin{array}{l}
{\displaystyle
\frac{1}{2}\binom{n+2}{3} - \frac{1}{4}\binom{n+1}{2}
- \frac{1}{8}\binom{n}{1} - \frac{1}{8}
}
\hspace{2ex}\mbox{if $n$ is even,} \\
   \\
{\displaystyle
\frac{1}{2}\binom{n+2}{3} - \frac{1}{4}\binom{n+1}{2}
- \frac{1}{8}\binom{n}{1}
}
\hspace{2ex}\mbox{if $n$ is odd.}
\end{array}\right.
\]
\end{thm}

\end{document}